\documentclass[11pt]{article}

\usepackage{amsmath}
\usepackage{amssymb}
\usepackage{amsthm}
\usepackage{latexsym}
\usepackage{color}
\usepackage{graphicx}
\usepackage{color}
\usepackage{caption}
\usepackage{subcaption}
\usepackage{hyperref}
\usepackage{enumerate}
\usepackage{mathtools}
\usepackage{accents}
\usepackage{ulem}

\DeclareSymbolFont{calletters}{OMS}{cmsy}{m}{n}
\DeclareSymbolFontAlphabet{\mathcal}{calletters}

\newtheorem{Theorem}{Theorem}[section]
\newtheorem{Definition}[Theorem]{Definition}
\newtheorem{Proposition}[Theorem]{Proposition}

\newtheorem{Lemma}[Theorem]{Lemma}

\newtheorem{Remark}[Theorem]{Remark}

\makeatletter \@addtoreset{equation}{section}

\newcommand{\brak}[1]{\left(#1\right)}    % round brackets
\newcommand{\crl}[1]{\left\{#1\right\}}   % curly brackets
\newcommand{\edg}[1]{\left[#1\right]}     % edgy brackets
\newcommand{\abv}[1]{\left|#1\right|}     % absolute value
\newcommand{\norm}[1]{\left\|#1\right\|}  % norm
\newcommand{\agb}[1]{\left\langle #1\right\rangle}   %qv

\usepackage{color}

\newcommand{\bm}{\bibitem}
\newcommand{\no}{\noindent}

\def \E{\mathbb{E}}

\def \R{\mathbb{R}}
\def \S{\mathbb{S}}

\def \N{\mathbb{N}}

\def\Ac{{\cal A}}
\def\Bc{{\cal B}}

\def\Hc{{\cal H}}

\def\Nc{{\cal N}}

\def\Pc{{\cal P}}

\def\Sc{\mathcal{S}}

\def\Uc{{\cal U}}

\def\Cb{\overline{C}}

\def\Et{\widetilde{\E}}

\def\Yt{\tilde{Y}}
\def\Zt{\tilde{Z}}
\def\Lv{\vec{L}}

\def\Ebh{\widehat{\E}}

\def \sigmal{\underline{\sigma}}
\def \a{\alpha}
\def \ah{\frac{\alpha}{2}}

\def \b{\beta}

\def \Om{\Omega}
\def \om{\omega}

\def \xb{\mathbf{x}}

\def \0{\mathbf{0}}

\def \Yb{\overline{Y}}
\def \Yh{\hat{Y}}

\def \Zh{\hat{Z}}
\def \Zb{\overline{Z}}
\def \Kb{\overline{K}}

\def \Kt{\tilde{K}}

\def \x{\times}
\def \sigmah{\overline{\sigma}}

\def \Dh{\delta}
\def \GH{\Gamma}

\def \e{\mathbf{e}}
\def \1{\mathbf{1}}
\def \xb{\mathbf{x}}
\def \Lm{\Lambda}
\def \Lt{\tilde{L}}
\def \Ht{\tilde{H}}

\usepackage{authblk}
\title{Mean-reflected $G$-BSDEs with multi-variate constraints}
\author[a, b]{Yiqing Lin}
\author[c]{Falei Wang}
\author[b]{Hui Zhao\thanks{Corresponding author: Hui Zhao (amyzh9201@sjtu.edu.cn)}}
\affil[a]{MOE-LSC, Shanghai Jiao Tong University, Shanghai, 200240, China}
\affil[b]{School of Mathematical Sciences, Shanghai Jiao Tong University, Shanghai, 200240, China}
\affil[c]{Zhongtai Securities Institute for Financial Studies, Shandong University, Jinan, 250100, China }

\begin{document}

\maketitle

\begin{abstract}
    In this paper, we study the multi-dimensional reflected backward stochastic differential equation driven by $G$-Brownian motion ($G$-BSDE) with a multi-variate constraint on the $G$-expectation of its solution. The generators are diagonally dependent on $Z$  and on all $Y$-components. We obtain the existence and uniqueness result via a fixed-point argumentation.

	\vspace{2mm}

	\no {\bf Key words.} Nonlinear expectation; $G$-Brownian motion; Backward stochastic differential equations; Mean reflection.

\end{abstract}

 \section{Introduction}
In 1990s, Pardoux and Peng \cite{PengBSDE} established the theory of nonlinear backward stochastic differential equations (BSDEs) of the form
$$Y_t=\xi+\int_t^Tf(s,Y_s,Z_s)ds-\int_t^TZ_sdB_s,$$
to which the existence and uniqueness of an adapted solution $(Y,Z)$ is obtained when the generator $f$ is uniformly Lipschitz and the terminal value $\xi$ is square-integrable. Since then researchers have developed the BSDE theory extensively and application of BSDEs are widely found in the domains of partial differential equations (PDEs), stochastic controls and mathematical finance, see  \cite{bsdefinanc, bsdepde, Pengbsdecontrol} for more details.\\

Among those extensions, reflected BSDEs (RBSDEs) associated with obstacles and the accompanied optimal stopping problems was introduced  by El Karoui, Kapoudjian, Pardoux, Peng and Quenez in \cite{rbsde}. Typically, the equations are of the following form
\begin{align*}
	Y_t=\xi+\int_t^Tf(s,Y_s,Z_s)ds-\int_t^TZ_sdB_s+K_T-K_t,\ Y_t\geq L_t, \ 0\leq t\leq T,
\end{align*}
where $L$ is a given obstacle process and $K$ is a non-decreasing process giving compensation to the equation in a minimal way that it satisfies the Skorokhod condition i.e., $\int_0^T\brak{Y_t-L_t}dK_t=0$. For more information, one may refer to \cite{convexbsde,SupL, PenRB, 2rBSDE} and the references therein.
Apart from pathwise constraints on the solution, BSDEs with constraints on distributions have been attracting a great deal of attentions. Bouchard, Elie and Réveillac \cite{bsdewt} introduced a type of BSDE with weak terminal condition, that is, instead of a concrete random variable, the terminal condition is imposed with distributions. Inspired from the risk management problems, Briand, Elie and Hu \cite{mrbsdeBEH} introduced the following mean-reflected BSDEs:
$$Y_t=\xi+\int_t^Tf(s,Y_s,Z_s)ds-\int_t^TZ_sdB_s+R_T-R_t,~ \E\edg{\ell(t,Y_t)}\geq 0, ~ 0\leq t\leq T,$$
in which $\ell(t,\cdot)$ is a family of running loss functions. According to \cite{mrbsdeBEH}, $(Y, Z, R)$ is called a flat solution of the above equation when $R$ is a deterministic process verifying the Skorokhod condition $\int_0^T\E\edg{\ell(t,Y_t)}dR_t=0$. In unidimensional setting, the well-posedness is proved with Lipschitz drivers and square-integrable terminal values in \cite{mrbsdeBEH}. Thereafter, the research of mean-reflected SDEs and BSDEs are developed in various directions, for instance \cite{ParticleMR,MRQBSDE, GeMRBSDE} etc..

 Motivated by financial problems with Knightian uncertainty, Peng \cite{PengMdG, PengGIto,PengNEC} has constructed systematically a framework of time-consistent sublinear expectation, named $G$-expectation theory, in which a new type of Brownian motion associated with the nonlinear expectation is introduced. Briefly speaking, the $G$-Brownian motion is a non-trivial generalization of its counterpart in the classical framework, so is the related stochastic calculus,  see  e.g. \cite{HuJiLiu,HuPeng,LPSSuperD,LiXPeng, martG,SongGeva}. Moreover, Hu, Ji, Peng and Song \cite{HuJiPengSong} demonstrated the well-posedness of one-dimensional backward stochastic differential equations driven by $G$-Brownian motion ($G$-BSDEs):
$$Y_t=\xi+\int_t^Tf(s,Y_s,Z_s)ds+\int_t^Tg_{ij}(s,Y_s,Z_s)d\agb{B^i,B^j}_s-\int_t^TZ_sdB_s-(K_T-K_t),$$
where $K$ is a non-increasing $G$-martingale.
Furthermore, one can turn to \cite{ComFKGir}, the sequel of \cite{HuJiPengSong}, for the further study of the comparison theorem, the Feynman-Kac formula and the Girsanov transformation. As for the multi-dimensional case, Liu \cite{LiumdGBSDE} obtained the solvability result for $G$-BSDEs with diagonal generators.

On reflected $G$-BSDEs, Li, Peng and Soumana Hima \cite{1dRGBSDE} investigated the one-dimensional problem with lower obstacles: 
$$Y_t=\xi+\int_t^Tf(s,Y_s,Z_s)ds+\int_t^Tg(s,Y_s,Z_s)d\agb{B}_s-\int_t^TZ_sdB_s+(A_T-A_t), \quad Y_t\geq S_t, \quad t\in [0,T],$$
where $A$ is a non-decreasing process and their flatness is given by an alternative criterion that the process  $-\int_0^\cdot\brak{Y_s-S_s}dA_s $ is a non-increasing $G$-martingale. Thereafter, Li and Liu \cite{liLiuRGd} developed the results for multi-dimensional reflected $G$-BSDEs. Moreover, the well-posedness of reflected $G$-SDEs and $G$-BSDEs can also be found in the works \cite{rqBSDE,UpOb, LinGSDE}.

In addition to pathwise obstacle problems in the $G$-framework, mean-reflected $G$-BSDEs of the following structure have been studied in \cite{LWMFGB, GLX}:
\begin{equation*}
	\begin{cases}
		Y_t=\xi+\int_t^Tf(s,Y_s,Z_s)ds+\int_t^Tg(s,Y_s,Z_s)d\agb{B}_s-\int_t^TZ_sdB_s-(K_T-K_t)+(R_T-R_t);\\
		\Ebh\edg{\ell(t,Y_t)}\geq 0, ~ 0\leq t\leq T,
	\end{cases}
\end{equation*}
where $\Ebh\edg{\cdot}$ is a nonlinear expectation measuring the uncertainty of running loss and a Skorokhod condition is given by $\int_0^T\Ebh\edg{\ell(t,Y_t)}dR_t=0$. 
Liu and Wang \cite{LWMFGB} studied the existence and uniqueness of the solution $(Y,Z,K,R)$ with Lipschitz generators and bounded terminal values. Gu, Lin and Xu \cite{GLX} worked on quadratic mean-reflected $G$-BSDEs when $Z$ is in the space $BMO_G$ and established solvability results for unbounded terminal conditions. For the regards on ``worst-case scenario", the constraint in \cite{GLX} was modified as $\Ebh\edg{-\ell(t,Y_t)}\leq 0$ and the flatness was given as $\int_0^T\Ebh\edg{-\ell(t,Y_t)}dR_t=0$. Recently, W. He \cite{mdGnonlip} investigated multi-dimensional mean-reflected $G$-BSDEs with constraints related to $G$-expectation however, individually imposed for each component respectively.
	
	 It is worth mentioning that our paper is a first attempt on the case where the reflecting boundary involves interactions among all solution components.
In our setting, multi-dimensional mean-reflected $G$-BSDEs with constraint on nonlinear expectation vector is studied, where $\Ebh\edg{Y}$ is restricted in the half space $\crl{y\in\R^N: \sum_{l=1}^N\theta^l y^l\geq 0}$ with $\theta^l ,~1\leq l\leq N$, being coefficients of a convex combination: 
\begin{align*}
		Y_t^l=&\ \xi^l +\int_t^T f^l (s,Y_s,Z_s^l)ds+\int_t^T g_{i,j}^l (s,Y_s,Z_s^l)d\agb{B^{i},B^{j}}_s\\
		&-\int_t^T Z_s^ldB_s-\brak{K_T^l-K_t^l}+\brak{R_T^l-R_t^l},\ 1\leq l\leq N.
\end{align*}
In view of ``worst-case scenario", if $\sum_{l=1}^N\theta^l \Ebh\edg{-Y_t^l}\leq 0$, then it holds that $\sum_{l=1}^N\theta^l \Ebh\edg{Y^l_t}\geq 0$. Accordingly, the flatness condition is given by $\int_0^T \brak{\sum_{l=1}^N\theta^l \Ebh\edg{-Y_t^l}}d\abv{R_t}=0$. Assumptions of Lipschitz growth and diagonal form on generators and of boundedness on terminal values are required for the fixed-point argument.     

This paper is organized as follows. In Section 2, we recall some preliminary knowledge of the $G$-framework, in particular, the corresponding $G$-BSDE theory. In Section 3, we construct the multi-variate mean-reflected $G$-BSDEs. In Section 4, we present our well-posedness results. And in Section 5, we explore an extension to reflection problems with constraints on general nonlinear expectations.

\section{Preliminaries}

\subsection{$G$-expectation framework}
We start by reviewing some essential notations and theories under the $G$-framework initiated by  Peng \cite{PengNEC}.
Let $\Om$ be a complete, separable metric space and let  $\Hc$ be a linear space of real-valued functions defined on $\Om$, such that $c\in\Hc,\ \forall c\in\R$ and $\abv{X}\in\Hc,\ \forall X\in\Hc$.
\begin{Definition}[\cite{PengNEC}.\ Sublinear expectation]
	A sublinear expectation $\Ebh:\Hc\to\R$ is a functional satisfying: for all $X,\ Y\in \Hc$,
	\begin{enumerate}[(1)]
		\item \textbf{Monotonicity:} $\Ebh\edg{X}\geq\Ebh\edg{Y}$, if $X\geq Y$;
		\item \textbf{Sub-additivity:} $\Ebh\edg{X+Y}\leq \Ebh\edg{X}+\Ebh\edg{Y}$;
		\item \textbf{Positive homogeneity:} $\Ebh\edg{\lambda X}=\lambda\Ebh\edg{X}$, for all $\lambda\geq 0$;
		\item \textbf{Constant preservation:} $\Ebh\edg{X+c}=\Ebh\edg{X}+c$, for all $c\in\R$.
	\end{enumerate} 
The triple $(\Om,\Hc,\Ebh)$ is called a sublinear expectation space.
\end{Definition} 

\begin{Definition}[\cite{PengNEC}. Independence]
	In a sublinear expectation space $(\Om,\Hc, \Ebh)$, a random vector $Y=\left(Y_1,\right.$$\cdots$,$\left.Y_n\right)^T$ $\in\Hc^n$ is said to be independent from another random vector $X=\brak{X_1,\cdots,X_m}^T\in\Hc^m$ under $\Ebh\edg{\cdot}$, if for each test function $\varphi\in C_{l,Lip}(\R^{m+n})$ we have
$$\Ebh\edg{\varphi(X,Y)}=\Ebh\edg{\Ebh\edg{\varphi(x,Y)}_{x=X}},$$
where $C_{l,Lip}(\R^n)$ is the linear space of continuous functions defined on $\R^n$ such that
$$\abv{\varphi(x)-\varphi(y)}\leq C\brak{1+\abv{x}^k+\abv{y}^k}\abv{x-y},$$
 for $x$, $y\in\R^n$ and some $C>0$, $k\in\N$ depending only on $\varphi$.
\end{Definition}

\begin{Definition}[\cite{PengNEC}. $G$-normal distribution $\Nc(\{0\}\x\Sigma)$]
	
	The random vector $X=\brak{X_1,X_2,\cdots, X_d}^T$ is said to be $G$-normally distributed, if for any $\varphi\in C_{l,Lip}(\R^d)$, the function $u(t,x):=\Ebh\edg{\varphi(x+\sqrt{t}X)}$ is a viscosity solution to the $G$-heat equation on $\R^+\x\R^d$:
$$\partial_tu-G(\partial_{xx}^2 u)=0,\ u(0,x)=\varphi(x).$$
Here $G:\S_d\to\R$ is a monotonic, sublinear mapping defined by
	$G(A):=\frac12\Ebh\edg{\agb{AX,X}}$, where $\S_d$ denoting the space of $d\x d$ symmetric matrices.
\end{Definition}
By Peng \cite{PengNEC}, there exists a bounded, convex and closed subset $\Sigma\subset\S_d^+$ such that
$G(A)=\frac12\sup_{X\in\Sigma}\text{tr}(XA),$
where $\S_d^+$ denotes the non-negative definite elements of $\S_d$. 
In this paper, we will always assume that $G$ is non-degenerate, i.e. $G(A)-G(B)\geq \sigmal^2\text{tr}\edg{A-B}$ for some $\sigmal>0$.

The $G$-expectation space $(\Om_T, L_{ip}(\Om_T),\Ebh,(\Ebh_t)_{t\geq 0})$ for some deterministic terminal time $T>0$ is constructed as follows.
 Let $\Om_T:= C_0 ([0,T], \R^d)$ be the space of continuous paths on $[0,T]$ starting from the origin, equipped with the supremum norm, $\left\|\om\right\|:=\sup_{t\in [0,T]}\abv{\om_t}$. Set $\Om_t:=\{\om_{\cdot\wedge t}: \om\in\Om\}$ and denote $\Bc(\Om_T)$ (resp. $\Bc(\Om_t)$) the Borel $\sigma$-algebra of $\Om_T$ (resp. $\Om_t$). For each $t\in [0,T]$, define
$$L_{ip}(\Om_t):=\crl{\varphi(B_{t_1},...,B_{t_k}):k\in\N,\ 0\leq t_1<\cdots<t_k\leq t,\ \varphi\in C_{b,Lip}(\R^{k\x d})},$$
where $C_{b,Lip}(\R^{k\times d})$ denotes the space of bounded and Lipschitz functions on $\R^{k\x d}$. 

\begin{Definition}[\cite{PengNEC}. $G$-Brownian motion and $G$-expectation]
	A $d$-dimensional process $\brak{B_t}_{t\geq 0}$ on a sublinear expectation space $(\Om, L_{ip}(\Om_T),\Ebh)$ is called a $G$-Brownian motion if 
	\begin{enumerate}[(i)]
		\item $B_0(\om)=\0$;
		\item for each $0<t\leq s\leq T$, the increment $B_{s}-B_{t}$ is $\Nc(\{0\}, (s-t)\Sigma)$-distributed and is independent from $(B_{t_1},B_{t_2},\cdots,B_{t_k})$, for any $k\in\N$ and $0\leq t_1\leq\cdots\leq t_k\leq t$.
	\end{enumerate}
    The sublinear expectation $\Ebh\edg{\cdot}$ is called a $G$-expectation if the canonical process $B$ is a $G$-Brownian motion under $\Ebh\edg{\cdot}$.
\end{Definition}

For each $\xi\in L_{ip}(\Om_T)$ taking the form
$$\phi\brak{B_{t_1}, B_{t_2}-B_{t_1},\cdots,B_{t_n}-B_{t_{n-1}}},\ \phi\in C_{b,Lip}(\R^{n\x d}),$$
the conditional $G$-expectation $\Ebh_{t_i}\edg{\cdot}, i=1,2,\cdots,n$, is defined in the following way:

\centerline{$\Ebh_{t_i}\edg{\phi\brak{B_{t_1}, B_{t_2}-B_{t_1},\cdots,B_{t_n}-B_{t_{n-1}}}}=\psi\brak{B_{t_1}, B_{t_2}-B_{t_1},\cdots,B_{t_i}-B_{t_{i-1}}},$}

\noindent where 
$\psi(x_1,x_2,\cdots,x_i)=\Ebh\edg{\phi(x_1,x_2,\cdots,x_i,B_{t_{i+1}}-B_{t_i},\cdots,B_{t_n}-B_{t_{n-1}})}.$
If $t\in(t_i,t_{i+1})$, we can reformulate $\xi$ as

\centerline{$\xi=\hat{\phi}\brak{B_{t_1}, B_{t_2}-B_{t_1},\cdots,B_{t}-B_{t_{i}},B_{t_{i+1}}-B_t,\cdots, B_{t_n}-B_{t_{n-1}}},\ \hat{\phi}\in C_{b,Lip}(\R^{d\x(n+1)}).$}

For $p\geq 1$, denote the norm $\norm{\xi}_{L_G^p}:=\Ebh\edg{\abv{\xi}^p}^{\frac1p}$ and $L_G^p(\Om_T)$ the completion of $L_{ip}(\Om_T)$ under $\norm{\cdot}_{L_G^p}$. Then the $G$-expectation $\Ebh\edg{\cdot}$ and the conditional $G$-expectation $\Ebh_t\edg{\cdot}$ can be continuously extended to $L_G^p(\Om_T)$.

\begin{Definition}[\cite{PengNEC}. $G$-martingale]
	A process $\{M_t\}_{t\in [0,T]}$ is called a $G$-martingale if
	\begin{enumerate}[(i)]
		\item $M_t\in L_G^1(\Om_t)$, for any $t\in [0,T]$;
		\item $\Ebh_s\edg{M_t}=M_s$, for all $0\leq s\leq t\leq T$.
	\end{enumerate}
    The process $\{M_t\}_{t\in [0,T]}$ is called a symmetric $G$-martingale if $-M$ is also a $G$-martingale.
\end{Definition}

\begin{Theorem}[\cite{DenisHuPeng, HuPeng,PengNEC}. Representative set of $G$-expectation]
	There exists a weakly compact set $\Pc$ of probability measures on $(\Om_T,\Bc(\Om_T))$ such that
	$$\Ebh\edg{X}=\sup_{P\in\Pc}\E^P\edg{X},~\forall X\in L_G^1(\Om_T),$$
	where $\E^P\edg{\cdot}$ is the expectation under probability $P$. Such $\Pc$ is called a set that represents $\Ebh$.
\end{Theorem}    
Then the associated capacity $c(\cdot)$ is defined by
 $   c(A):=\sup_{P\in\Pc}P(A),~ A\in\Bc(\Om_T).$
A set $A\in\Bc(\Om_T)$ is a polar set if $c(A)=0$ and a property holds "quasi-surely" (q.s.) if it holds outside a polar set.
We do NOT distinguish the random variables $X$ and $Y$ satisfying $X=Y$ q.s..

\begin{Definition}
    Given $p\geq 1$, the following spaces will be frequently used in our framework:
    \begin{enumerate}[(i)]
    	\item $M_G^0(0,T):=\crl{\eta : \eta_t(\om)=\sum_{j=0}^{n-1}\xi_j(\om)1_{\left [t_j,t_{j+1}\right )}(t),~\mbox{where}~0=t_0<\cdots< t_n=T~\mbox{and}~\xi_j\in L_{ip}(\Om_{t_j})}$;
    	\item $M_G^p(0,T)$: the completion of $M_G^0(0,T)$
    	under the norm $\norm{\eta}_{M_G^p(0,T)}:=\crl{\Ebh\edg{\int_0^T\abv{\eta_s}^p ds}}^{\frac1p}$;
    	\item $H_G^p(0,T)$: the completion of $M_G^0(0,T)$ under the norm $\norm{\eta}_{H_G^p(0,T)}:=\crl{\Ebh\edg{\brak{\int_0^T\abv{\eta_s}^2 ds}^{p/2}}}^{\frac1p}$;
    	\item $S_G^0(0,T):=\crl{h(t,B_{t_1\wedge t},...,B_{t_n\wedge t}):t_1,...,t_n\in[0,T],\ h\in C_{b,Lip}(\R^{1+d\x n})}$;
    	\item $S_G^p(0,T)$: the completion of $S_G^0(0,T)$ under the norm $\norm{\eta}_{S_G^p(0,T)}:=\crl{\Ebh\edg{\sup_{t\in[0,T]}\abv{\eta_t}^p}}^{\frac1p}$;
    	\item $\Ac_G^p(0,T)$: the collection of processes $K$ such that $K$ is a non-increasing $G$-martingale with $K_0 =0$ and $K_T\in L_G^p(\Om_T)$;
    	\item $\Ac_D(0,T)$: the collection of deterministic, continuous and non-decreasing processes starting from 0;
    	\item $
    	\Sc_G^p(0,T):= S_G^p(0,T)\x H_G^p(0,T;\R^d)\x \Ac_G^p(0,T).
    	$	
    \end{enumerate}
\end{Definition}    
    
   For multi-dimensional case, we denote by $M_G^p(0,T;\R^N)$,
   the space of $N$-dimensional processes $\crl{\eta_t}_{t\in[0,T]}$ such that $\eta^l\in M_G^p(0,T),~\forall 1\leq l\leq N$.
   Similarly we introduce the notations as $L_G^p(\Om_T;\R^N)$, 
   $S_G^p(0,T;\R^N)$, $H_G^p(0,T;\R^{N\times d})$, $\Ac_G^p(0,T;\R^N)$, $\Ac_D(0,T;\R^N)$, and $\Sc_G^p(0,T;\R^N)$ for the corresponding spaces of $N$-dimensional random vectors and processes on $[0,T]$. 
   
   Given some time interval $[t_1,t_2]\subset[0,T]$, denote 
   $\Sc_G^\a(t_1,t_2;\R^N)=S_G^\a(t_1,t_2;\R^N)\x H_G^\a(t_1,t_2;\R^{N\x d})\x\Ac_G^\a(t_1,t_2;\R^N).$
   Let $\Ac_D(t_1,t_2;\R^N)$ be the space of $N$-dimensional deterministic, continuous and non-decreasing processes being $\0$ at $t_1$.

    \begin{Lemma}[\cite{DenisHuPeng}. Monotone convergence and Fatou's lemma]\label{lemma:FatouMonotone}
    	Let $\crl{X_n}_{n\geq 1}$ be a sequence of $\Bc(\Om_T)$-measurable random variables.
    	\begin{enumerate}
    		\item[(i)] Suppose $X_n\geq 0$. Then $\Ebh\edg{\liminf_{n\to\infty}X_n}\leq\liminf_{n\to\infty}\Ebh\edg{X_n}$.
    		\item[(ii)] Suppose $\brak{X_n}_{n\geq 1}\subset L_G^1(\Om_T)$ are non-increasing. Then $\Ebh\edg{X_n} \downarrow\Ebh\edg{\lim_{n\to\infty}X_n}$.
    	\end{enumerate}
    \end{Lemma}

Note that the $\limsup$ version of Fatou's lemma does NOT hold in the $G$-framework, as a consequence, the dominated convergence theorem does NOT hold in general with the norm $\norm{\cdot}_{L_G^1}$.  Given the $G$-Brownian motion $B$, denote by $\agb{B^i,B^j},~1\leq i,j\leq d$ the mutual variation processes. Then for  $\eta\in H_G^p(0,T)$ and $\zeta\in M_G^p(0,T)$ with $p\geq 1$, the $G$-It\^o integrals $\int_0^\cdot\eta_sdB^i_s$ and $\int_0^\cdot \zeta_sd\agb{B^i,B^j}_s$ are well-defined. 
    \begin{Lemma}[\cite{SongGeva}. BDG type inequality]\label{lemma:BDG}
    	For $\eta\in H_G^\a(0,T)$ with $\a\geq 1$ and $p\in\left(0,\a\right]$, we have 
    	
    	$$\sigmal^pc_p\Ebh\edg{\brak{\int_0^T\abv{\eta_r}^2dr}^{\frac{p}2}}\leq \Ebh\edg{\sup_{t\in[0,T]}\abv{\int_0^t\eta_rdB_r}^p}\leq \sigmah^pC_p\Ebh\edg{\brak{\int_0^T\abv{\eta_r}^2dr}^{\frac{p}2}},~\mbox{q.s.,}$$ 
    	where $c_p,C_p>0$ are independent of $\eta$.
    \end{Lemma}
     
     \begin{Theorem}[\cite{PengNEC}. Doob's maximal inequality for $G$-martingales]\label{theorem:simpleDoob}
     	Suppose $1\leq\alpha<\beta$. For each $1<p<\overline{p}:=\beta/\alpha$ with $p\leq 2$ and for all $\eta\in L_G^\beta(\Om_T)$, there exists a constant $C>0$ depending only on $p$, $\overline{\sigma}$ and $\underline{\sigma}$ such that
     	$$\Ebh\edg{\sup_{t\in[0,T]}\Ebh_t\edg{\abv{\eta}^\alpha}}\leq\frac{Cp\overline{p}}{\brak{\overline{p}-p}\brak{p-1}}\brak{\Ebh\edg{\abv{\eta}^\beta}^{\frac{1}{p\overline{p}}}+\Ebh\edg{\abv{\eta}^\beta}^{\frac{1}{p}}}.$$
     \end{Theorem}

\begin{Lemma}[\cite{LWMFGB}. Uniform estimate for conditional $G$-expectation]\label{lemma:SGui}
	Suppose that $X\in L_G^p(\Om_T)$ with $p\geq 1$. Then, for any $\varepsilon>0$, there exists a constant $\delta>0$ such that for all $A\in\Bc(\Om_T)$ and $c(A)\leq\delta$, we have
	$$\sup_{t\geq 0}\Ebh\edg{\abv{\Ebh_t\edg{X}}^p\1_A}\leq \varepsilon.$$
\end{Lemma}
\begin{Lemma}[\cite{HuJiLiu}. Continuity property of conditional $G$-expectation on $L_{ip}(\Om_T)$]\label{lemma:CELipcon}
Let $X\in L_{ip}(\Om_T)$, then for any $T\geq 0$ and $0\leq s_1\leq s_2\leq T$, we have
\begin{equation*}
	\abv{\Ebh_{s_2}\edg{X}-\Ebh_{s_1}\edg{X}}\leq C\Big\{\sup_{(u_1,u_2)\in\Lm_{s_2-s_1,T}}\brak{\abv{B_{u_2}-B_{u_1}}\wedge1}+\sqrt{s_2-s_1}\Big\},
\end{equation*} 	
where $C$ is a constant depending only on $X$ and $G$, and
$$	\Lm_{\delta,T}:=\crl{(u_1,u_2):0\leq u_1,u_2\leq T,~\abv{u_1-u_2}\leq \delta}.$$
Moreover, 
$$\Ebh\Big[\sup_{(u_1,u_2)\in\Lm_{\delta,T}}\brak{\abv{B_{u_2}-B_{u_1}}\wedge1}\Big]\downarrow0,~\mbox{as}~~\delta\downarrow0.$$	
\end{Lemma}	
It is evident that the multi-dimensional version of Lemma \ref{lemma:CELipcon} holds as well.
\begin{Lemma}[\cite{LPSSuperD}. Continuity property of $S_G^p(0,T)$]\label{lemma:SGcon}
	For $Y\in S_G^p(0,T)$ with $p\geq 1$, we have
	$$F(Y):=\limsup_{\varepsilon\to 0}\brak{\Ebh\edg{\sup_{t\in[0,T]}\sup_{s\in\edg{t,(t+\varepsilon)\wedge{T}}}\abv{Y_t-Y_s}^p}}^{\frac{1}{p}}=0.$$
\end{Lemma} 	
\begin{Lemma}[\cite{LiumdGBSDE}. Characterization of $M_G^p(t,t+h)$]\label{lemma:funcappr}
	Suppose $p\geq 1$. Given a function $f:[t,t+h]\x\Om_T\x\R^m\to\R$ satisfying $f(\cdot,\cdot,x)\in M_G^p(t,t+h)$ for each $x\in \R^m$. Assume that $f$ is Lipschitz continuous in $x$ with some coefficient $L>0$, uniformly in $(s,\om)$. Then for any $X\in M_G^p(t,t+h;\R^m)$, we have $(f(s, \cdot, X_s(\cdot))_{t\leq s\leq t+h}\in M_G^p(t,t+h)$. 
\end{Lemma}

\subsection{$G$-BSDE theory}

 Next we shall recall some fundamental results on one-dimensional $G$-BSDEs (cf. Hu, Ji, Peng and Song \cite{HuJiPengSong}) of the type (Einstein's convention is applied):
\begin{equation}\label{eq:1dGBSDE}
	Y_t=\xi+\int_t^Tf(s,Y_s,Z_s)ds+\int_t^Tg_{ij}(s,Y_s,Z_s)d\agb{B^i,B^j}_s-\int_t^TZ_sdB_s-(K_T-K_t),
\end{equation}
where the generators
$f$, $g_{ij}: [0,T]\x\Om_T\x\R\x\R^d\to\R,$
 satisfy:
\begin{enumerate}
	\item[($A_1$)] there exists some $\beta>1$ such that for any $(y,z)\in\R\x\R^d$, $f(\cdot,\cdot,y,z),g(\cdot,\cdot,y,z)\in M_G^\beta(0,T)$;
	\item[($A_2$)] there exists some $L>0$ such that for any $y,y'\in\R,z,z'\in\R^d$, 
	
	$\abv{f(t,\om,y,z)-f(t,\om,y',z')}+\sum_{i,j=1}^d\abv{g_{i,j}(t,\om,y,z)-g_{ij}(t,\om,y',z')}\leq L\brak{\abv{y-y'}+\abv{z-z'}};$
	\item[($A_\xi$)] $\xi\in L_G^\beta(\Om_T)$.	
\end{enumerate}

\begin{Theorem}[\cite{HuJiPengSong}]\label{thm:1deu}
	Assume that $\xi, f$ and $g_{i,j}$ satisfy ($A_1$)-($A_2$)-($A_\xi$) for some $\beta>1$ and $L>0$. Then for any $1<\a<\b$, $G$-BSDE \eqref{eq:1dGBSDE} admits a unique solution $(Y,Z,K)\in \Sc_G^\a(0,T)$.
\end{Theorem}

Here are some a priori estimates.

\begin{Proposition}[\cite{HuJiPengSong}]\label{prop:1daprY}
	Let $\xi$, $f$ and $g_{i,j}$ satisfy ($A_1$)-($A_2$)-($A_\xi$) for some $\beta>1$ and $L>0$. Assume that $(Y,Z,K)\in \Sc_G^\a(0,T)$, for some $1<\a\leq\b$, is a solution to equation \eqref{eq:1dGBSDE} with parameter $(\xi,f,g_{ij})$. Then there exists some constant $C_{(\a)}>0$ depending on $T$, $G$, $L$ and $\a$ such that 
	$$\abv{Y_t}^\a\leq C_{(\a)}\Ebh_t\edg{\abv{\xi}^\a+\brak{\int_t^Th_sds}^\a},$$	
where $h_s=\abv{f(s,0,0)}+\sum_{i,j=1}^d\abv{g_{ij}(s,0,0)} $. 
\end{Proposition}

\begin{Proposition}[\cite{HuJiPengSong}]\label{prop:1daprDY}
	Let $\xi^{(m)}, f^{(m)}$ and $g_{i,j}^{(m)}, m=1,2$ satisfy ($A_1$)-($A_2$)-($A_\xi$) for the same $\beta>1$ and $L>0$. Assume that $(Y^{(m)},Z^{(m)},K^{(m)})\in \Sc_G^\a(0,T)$, for some $1<\a\leq\b$, is a solution to equation \eqref{eq:1dGBSDE} with parameter $(\xi^{(m)},f^{(m)},g_{ij}^{(m)})$. Set $\Yh_t:=Y^{(1)}_t-Y^{(2)}_t$. Then there exists some constant $C_{(\a)}>0$ depending on $T$, $G$, $L$ and $\a$ such that
	$$\abv{\Yh_t}^\a\leq C_{(\a)}\Ebh_t\edg{\abv{\hat{\xi}}^\a+\brak{\int_t^T\hat{h}_sds}^\a},$$	
where $\hat{\xi}=\xi^{(1)}-\xi^{(2)}$ and
$$\hat{h}_s=\abv{f^{(1)}(s,Y_s^{(2)},Z_s^{(2)})-f^{(2)}(s,Y_s^{(2)},Z_s^{(2)})}+\sum_{i,j=1}^d\abv{g_{ij}^{(1)}(s,Y_s^{(2)},Z_s^{(2)})-g_{ij}^{(2)}(s,Y_s^{(2)},Z_s^{(2)})}.$$
\end{Proposition}
\begin{Proposition}[\cite{HuJiPengSong}]\label{prop:1daprDZ}
	Let $\xi^{(m)}, f^{(m)}$ and $g_{i,j}^{(m)}, m=1,2$ satisfy ($A_1$)-($A_2$)-($A_\xi$) for the same $\beta>1$ and $L>0$. Assume that $(Y^{(m)},Z^{(m)},K^{(m)})\in \Sc_G^\a(0,T)$, for some $1<\a\leq\b$, is a solution to equation \eqref{eq:1dGBSDE} with parameter $(\xi^{(m)},f^{(m)},g_{ij}^{(m)})$. Set $\Yh_t:=Y^{(1)}_t-Y^{(2)}_t$ and $\Zh_t=Z^{(1)}_t-Z^{(2)}_t$. Then there exists some constant $C>0$ depending on $T$, $G$, $L$ and $\a$ such that 
	$$\Ebh\edg{\brak{\int_0^T\abv{\Zh_s}^2ds}^\ah}\leq C\crl{\norm{\Yh}_{S_G^\a}^\a+\norm{\Yh}_{S_G^\a}^\ah\sum_{m=1}^2\edg{\norm{Y^{(m)}}_{S_G^\a}^\ah+\norm{\int_0^Th_s^{(m)}ds}_{L_G^\a(\Om_T)}^\ah}},$$
	where $h_s^{(m)}=\abv{f^{(m)}(s,0,0)}+\sum_{i,j=1}^d\abv{g_{ij}^{(m)}(s,0,0)} $.
\end{Proposition}

\section{Multi-dimensional mean-reflected $G$-BSDE with multi-variate constraint}

\subsection{Formulation of the problem}

In this paper, by $x^l$ we denote the $l$-th component of $x\in\R^N$ ($\R^{N\x d}$) for $l=1,\cdots,N$. We will study the $N$-dimensional multi-variate mean-reflected $G$-BSDE (MR$G$-BSDE for short) taking the form:
\begin{equation}\label{eq:GBSDEDG}
	\left\{\begin{lgathered}
		Y_t^l=\xi^l +\int_t^T f^l (s,Y_s,Z_s^l)ds+\int_t^T g_{i,j}^l (s,Y_s,Z_s^l)d\agb{B^{i},B^{j}}_s-\int_t^T Z_s^ldB_s-\brak{K_T^l-K_t^l}\\
		\qquad~+\brak{R_T^l-R_t^l},	~1\leq l\leq N; 
		\\
	    \sum_{l=1}^N\theta^l \Ebh\edg{-Y_t^l}\leq 0,~0\leq t\leq T;~\int_0^T \brak{\sum_{l=1}^N\theta^l \Ebh\edg{-Y^l_t}} d\abv{R_t}=0.
	\end{lgathered}\right.	
\end{equation}
We assume the following assumptions:
\begin{enumerate}
	\item[($H_\theta$)]  $\theta^l \geq 0,~\forall 1\leq l\leq N$, and $\sum_{l=1}^N\theta^l =1$;
	\item[($H_\xi$)] the terminal value  $\xi\in L_G^\beta(\Om_T;\R^N)$ with $\b>2$ and $\sum_{l=1}^N\theta^l \Ebh\edg{-\xi^l}\leq 0$;

	\item[($H_f$)] the generators $ f^l(t,\om,y,z^l),~g_{ij}^l(t,\om,y,z^l):[0,T]\times\Om_T\times\R^N\times\R^d\to\R, ~ 1\leq l\leq N,$ satisfy that
	\begin{enumerate}[($i$)]
		\item for each $(y,z^l)\in\R^N\x\R^{d}$, $f^l(\cdot,\cdot,y,z^l),~ g_{ij}^l(\cdot,\cdot,y,z^l)\in M_G^\beta(0,T)$ with $\b>2$ in ($H_\xi$);
		\item there exists some $L>0$ such that for any $y_1,y_2\in\R^N$ and $z_1^l,z_2^l\in\R^d$,
		
		$\abv{f^l(t,y_1,z_1^l)-f^l(t,y_2,z_2^l)}+\sum_{i,j=1}^d \abv{g_{ij}^l(t,y_1,z_1^l)-g_{ij}^l(t,y_2,z_2^l)}\leq L\brak{\abv{y_1-y_2}+\abv{z_1^l-z_2^l}}.$
	\end{enumerate} 
\end{enumerate}

	We consider the constraint $\sum_{l=1}^N\theta^l \Ebh\edg{-Y_t^l}\leq 0$,  $0\leq t\leq T$, which is corresponding to a hedging problem with a constraint on the risk where $Y$ denotes the wealth. On the other hand, for the half-space $D:=\crl{y\in\R^N:\sum_{l=1}^N\theta^l y^l\geq 0}\subset\R^N$, an intuitive way to put the constraint on the mean vector  $\Ebh\edg{Y_.}:=(\Ebh\edg{Y^1_.},\cdots,\Ebh\edg{Y^N_.})^T$ is that $\Ebh\edg{Y_t}\in D,~ \forall t\in[0,T]$. We observe that, by properties of sublinear expectation,
	$$
	\sum_{l=1}^N\theta^l \Ebh\edg{Y_t^l}\geq\Ebh\edg{\sum_{l=1}^N\theta^l Y_t^l}\geq-\Ebh\edg{-\sum_{l=1}^N\theta^l Y_t^l}\geq-\sum_{l=1}^N\theta^l \Ebh\edg{-Y^l_t},~\forall t\in[0,T].
	$$
	Therefore, the constraint in \eqref{eq:GBSDEDG} 
	guarantees that $\Ebh\edg{Y_t}\in D, ~\forall t\in[0,T]$. However, there isn't any additional technique difficulty involved when either type of constraint is imposed.

For symbol simplicity in the rest of this paper, we note that $g(t,y,z):=\brak{g_{ij}(t,y,z)}_{i,j=1}^d$ and 
\centerline{$\psi(t,y,z):=\brak{\psi^1(t,y,z^1),\cdots,\psi^N(t,y,z^N)}^T$, for $\psi=f,\ g_{ij}$, $i$, $j = 1,\ldots d$.}

\begin{Theorem}\label{thm:existuniq}
	Suppose that ($H_\theta$)-($H_\xi$)-($H_f$) hold for some $\b>2$ and $L>0$. Then for any $2\leq\a<\beta$, MR$G$-BSDE \eqref{eq:GBSDEDG} with parameter $(\xi,f,g)$ admits a unique solution $(Y,Z,K,R)\in\Sc_G^\a(0,T;\R^N)\x\Ac_D(0,T;\R^N)$. Moreover, $Y\in M_G^\beta(0,T;\R^N)$.
\end{Theorem}

\subsection{Construction of the solution}

Similarly to the works by Liu and Wang \cite{LWMFGB} on one-dimensional MR$G$-BSDEs, we construct here a mapping $\Lv:L_G^1(\Om_T;\R^N)\to\R^N$:
\begin{align*}
	\Lv(X):=\arg\min\crl{\abv{\xb}:\xb\in \Uc(X)},~\mbox{with}~~\Uc(X):=\crl{\xb\in\R^N: \sum_{l=1}^N\theta^l \Ebh\edg{-X^l-\xb^l}\leq 0},
\end{align*} 
for $X\in L_G^1(\Om_T;\R^N)$. Then we denote by
$$H(\xb,X):=\sum_{l=1}^N\theta^l \Ebh\edg{-X^l-\xb^l}=-\sum_{l=1}^N\theta^l \xb^l+\sum_{l=1}^N\theta^l \Ebh\edg{-X^l}.$$
\begin{Remark}
	\begin{enumerate}[(i)]
		\item Note that the set  $~\Uc(X)\neq\emptyset$ for any $X\in L_G^1(\Om_T;\R^N)$, thus $\Lv(\cdot)$ is well-defined. 
		\item Obviously, $\Lv(\0)\equiv\0$.
	\end{enumerate}
\end{Remark}

Here we can obtain an explicit expression of $\Lv(X):=\brak{\Lv^1(X),\cdots,\Lv^N(X)}^T$ by solving the following quadratic programming problem with fixed $X$:
\begin{align*}
	&\mbox{minimize}~~ \frac12\abv{\xb}^2,\label{eq:QP}\\
	&\mbox{subject to}~~ \sum_{l=1}^N\theta^l  \xb^l=\brak{\sum_{j=1}^N\theta^j\Ebh\edg{-X^j}}^+.\notag
\end{align*}
Then we obtain the explicit formula: 
\begin{equation}\label{eq:L-l-t}
	\Lv^l(X):=\frac{\theta^l \brak{\sum_{j=1}^N\theta^j\Ebh\edg{-X^j}}^+}{\sum_{j=1}^N(\theta^j)^2},~1\leq l\leq N,~\mbox{and}~~\abv{\Lv(X)}=\frac{\brak{\sum_{j=1}^N\theta^j\Ebh\edg{-X^j}}^+}{\sqrt{\sum_{j=1}^N(\theta^j)^2}}.
\end{equation}

Hereafter, the following propositions will be frequently used.
\begin{Proposition}\label{prop:CEcont}
	Let $S\in L_G^1(\Om_T;\R^N)$ and $\zeta\in M_G^1(0,T;\R^N)$. Define  $\hat{S}:=\brak{\hat{S}^1,\cdots,\hat{S}^N}^T$ with $\hat{S}^l_t:=\Ebh_t\edg{S^l}-\int_0^t\zeta_s^lds,~ 1\leq l\leq N, ~t\in[0,T]$.
	\begin{enumerate}[(i)]
		\item The mapping $t\mapsto H(\xb,\hat{S}_t)$ is continuous.
		\item The mappings $t\mapsto \Lv(\hat{S}_t)$ and $t\mapsto\abv{\Lv(\hat{S}_t)}$ are continuous. 
	\end{enumerate}
\end{Proposition}
\noindent{\bf\underline{Proof.}}
$(i)$ Let $\xb$ be a fixed vector and $s<t$ w.l.o.g., we have
\begin{align*}
	\abv{H(\xb,\hat{S}_t)-H(\xb,\hat{S}_s)}
	&=\abv{\sum_{l=1}^N\theta^l \Ebh\edg{-\Ebh_t\edg{S^l}+\int_0^t\zeta_u^ldu}-\sum_{l=1}^N\theta^l \Ebh\edg{-\Ebh_s\edg{S^l}+\int_0^s\zeta_u^ldu}}\\
	&\leq \sum_{l=1}^N\theta^l \abv{\Ebh\edg{-\Ebh_t\edg{S^l}+\int_0^t\zeta^l_udu}-\Ebh\edg{-\Ebh_s\edg{S^l}+\int_0^s\zeta_u^ldu}}\\
	&\leq \sum_{l=1}^N\theta^l \Ebh\edg{\abv{\Ebh_t\edg{S^l}-\Ebh_s\edg{S^l}}}+\sum_{l=1}^N\theta^l \Ebh\edg{\int_0^T\1_{[s,t]}(u)\abv{\zeta_u^l}du}.
\end{align*}
Since $\zeta\in M_G^1(0,T;\R^N)$, we observe that the second term tends to 0 as $s\to t$. For the first term when $S\in L_{ip}(\Om_T;\R^N)$, we have by the multi-dimensional version of Lemma \ref{lemma:CELipcon}, 
\begin{equation*}
	\lim_{s\to t}\Ebh\edg{\abv{\Ebh_t[S^l]-\Ebh_s\edg{S^l}}}=0,\quad 1\leq l\leq N.
\end{equation*} 
For $S\in L_G^1(\Om_T;\R^N)$, the above convergence is trivial since $L_{ip}(\Om_T;\R^N)$ is dense in $L_G^1(\Om_T;\R^N)$.

$(ii)$ By the explict solution \eqref{eq:L-l-t}, one can find for $s<t$,
\begin{align*}
	\abv{\Lv^l(\hat{S}_t)-\Lv^l(\hat{S}_s)}\leq \frac{\theta^l }{\sum_{j=1}^N(\theta^j)^2}\sum_{j=1}^N\theta^j\crl{\Ebh\edg{\abv{\Ebh_t\edg{S^j}-\Ebh_s\edg{S^j}}}+\Ebh\edg{\int_0^T\1_{[s,t]}(u)\abv{\zeta^j_u}du}}.
\end{align*}
Following a similar procedure as ($i$), we can obtain the continuity result of $t\mapsto \Lv^l(\hat{S}_t)$.
\qed
 
\begin{Proposition}
	\begin{enumerate}[(i)]
		\item Let $S\in S_G^p(0,T;\R^N),~p\geq 1$. The mapping $t\mapsto \Lv(S_t)$ is continuous.
		\item For $X\in L_G^1(\Om_T;\R^N)$, the mapping $X\mapsto \Lv(X)$ is continuous under $\Ebh\edg{\abv{\cdot}}$.
	\end{enumerate}
\label{proposition:LconSG} 
\end{Proposition}
\noindent{\bf\underline{Proof.}}
$(i)$ For any $\varepsilon>0$, $0\leq s,t\leq T$ and each $1\leq l\leq N$, 
\begin{align*}
	\Ebh\edg{\abv{S^l_t-S^l_s}}&\leq \Ebh\edg{\abv{S^l_t-S^l_s}\1_{\crl{\abv{S^l_t-S^l_s}\leq\frac{\varepsilon}{2}}}}+\Ebh\edg{\abv{S^l_t-S^l_s}\1_{\crl{\abv{S^l_t-S^l_s}>\frac{\varepsilon}{2}}}}\\
	&\leq \frac{\varepsilon}{2}+2\Ebh\edg{\sup_{u\in[0,T]}\abv{S^l_u}\1_{\crl{\abv{S^l_t-S^l_s}>\frac{\varepsilon}{2}}}}.
\end{align*}
From Lemma \ref{lemma:SGcon} and the Markov's inequality, we have
$\lim_{s\to t}c\brak{\crl{\abv{S^l_t-S^l_s}>\frac{\varepsilon}{2}}}=0$. Thus, first for $S\in S^0_G(0, T)$, eventually for  $S\in S_G^p(0,T;\R^N)$, we can find a constant $\delta_l>0$ small enough and by Lemma \ref{lemma:SGui}, letting $\abv{t-s}<\delta_l$, to obtain that
$\Ebh\edg{\abv{S^l_t-S^l_s}}\leq\varepsilon.$

Then we choose $\delta:=\min_{1\leq l\leq N}\delta_l$ and put $\abv{s-t}<\delta$. From \eqref{eq:L-l-t}, it holds that
$$\abv{\Lv^l(S_t)-\Lv^l(S_s)}\leq \frac{\theta^l }{\sum_{j=1}^N(\theta^j)^2}\sum_{j=1}^N\theta^j\Ebh\edg{\abv{S^j_t-S^j_s}}\leq \frac{\theta^l \varepsilon}{\sum_{j=1}^N(\theta^j)^2}.$$
As $\varepsilon>0$ can be arbitrarily small, the proof ends.

$(ii)$ Given $X'\in L_G^1(\Om_T;\R^N)$, we have
\begin{align*}
	\abv{\Lv^l(X)-\Lv^l(X')}\leq\frac{\theta^l }{\sum_{j=1}^N(\theta^j)^2}\sum_{j=1}^N\theta^j\Ebh\edg{\abv{X^j-X^{'j}}}\leq \frac{\theta^l }{\sum_{j=1}^N(\theta^j)^2}\Ebh\edg{\abv{X-X^{'}}}.
\end{align*}
\qed

\section{Main results}
Before proving Theorem \ref{thm:existuniq}, we shall present two heuristic results concerning the well-posedness of MR$G$-BSDE in the following cases:
\begin{itemize}
	\item Case I: $N$-dimensional multi-variate MR$G$-BSDE with generators not depending on $(Y, Z)$;
	\item Case II: $N$-dimensional multi-variate MR$G$-BSDE with generators depending on $Z$ diagonally.
\end{itemize}
 In the sequel, we assume $g^l\equiv 0, ~1\leq l\leq N$ (the generator corresponding to the integral w.r.t. the co-variation of $G$-Brownian motion) and the results with general $g^l$ are similar.
     
\subsection{Case I}    
Consider an $N$-dimensional multi-variate MR$G$-BSDE with generators not depending on $(Y, Z)$. Then, \eqref{eq:GBSDEDG} can be rewritten as
\begin{equation}\label{eq:constdriver}
	\left\{\begin{lgathered}
		Y_t^l=\xi^l +\int_t^T f^l_sds-\int_t^T Z_s^l dB_s-\brak{K_T^l-K_t^l}+\brak{R_T^l-R_t^l},~\shoveright{1\leq l\leq N;} 
		\\
		\sum_{l=1}^N\theta^l \Ebh\edg{-Y_t^l}\leq 0, ~0\leq t\leq T;~\int_0^T \brak{\sum_{l=1}^N\theta^l \Ebh\edg{-Y^l_t}} d\abv{R_t}=0;
	\end{lgathered}	\right.
\end{equation}     
for which the assumption ($H_\theta$) remains and the assumptions ($H_\xi$)-($H_f$) are reduced to
\begin{itemize}
	\item[($H^I$)]  $\xi\in L_G^\beta(\Om_T;\R^N)$, $f\in M_G^\beta(0,T;\R^N)$ with $\b>2$, and $\sum_{l=1}^N\theta^l \Ebh\edg{-\xi^l}\leq 0$.
\end{itemize}    
\begin{Theorem}\label{theorem:constdriver}
	Suppose that ($H_\theta$)-($H^I$) holds. Then for any $~2\leq\a<\b$, MR$G$-BSDE \eqref{eq:constdriver} with parameter $(\xi,f)$ admits a unique solution $(Y,Z,K,R)\in\Sc_G^\a(0,T;\R^N)\x\Ac_D(0,T;\R^N)$. Furthermore, $Y\in M_G^\b(0,T;\R^N)$.
\end{Theorem}    

\noindent{\bf \underline{Proof.}} Upon proving flatness and uniqueness, the idea is borrowed from Proposition 7 in Briand, Elie and Hu \cite{mrbsdeBEH}.

\vspace{0.2em}
\noindent{\bf Step 1. Existence.}

For each $t\in[0,T]$, set $V_t=\brak{V_t^1,\cdots,V_t^N}^T$, where 
$$V^l_t:=\Ebh_t\edg{\xi^l+\int_t^Tf_s^lds}=\Ebh_t\edg{\xi^l+\int_0^Tf^l_sds}-\int_0^tf^l_sds,~1\leq l\leq N.$$  

From Proposition \ref{prop:CEcont}, both the mappings $t\mapsto H(\xb,V_t)$  and $t\mapsto \Lv(V_t)$ are continuous. Then we can define a deterministic continuous process $R\in\Ac_D(0,T;\R^N)$ by
\begin{equation}\label{eq:constR}
	R_t^l:=\sup_{s\in[0,T]}\Lv^l(V_s)-\sup_{s\in[t,T]}\Lv^l(V_s),~\mbox{and}~~\abv{R_t}:=\sup_{s\in[0,T]}\abv{\Lv(V_s)}-\sup_{s\in[t,T]}\abv{\Lv(V_s)}.
\end{equation}
For each $l\in\crl{1,\cdots,N}$, by Theorem \ref{thm:1deu}, the standard one-dimensional $G$-BSDE
\begin{equation}\label{eq;SDGBSDE}
	\Yb_t^l=\xi^l+\int_t^Tf^l_sds-\int_s^T\Zb_s^l dB_s-(\Kb_T^l-\Kb_t^l),~t\in[0,T],
\end{equation}   
admits a unique solution $(\Yb^l,\Zb^l,\Kb^l)\in\Sc_G^\a(0,T)$. Letting $Y_t^l=\Yb_t^l+(R_T^l-R_t^l)$, the quadruple of processes $$\brak{Y_t^l,Z_t^l,K_t^l,R_t^l}=\brak{\Yb^l_t+\sup_{s\in[t,T]}\Lv^l(V_t),~\Zb^l_t,~\Kb^l_t,\sup_{s\in[0,T]}\Lv^l(V_s)-\sup_{s\in[t,T]}\Lv^l(V_s)},$$
 is a solution to the one-dimensional reflected $G$-BSDE
\begin{equation}
	Y_t^l=\xi^l +\int_t^T f^l_sds-\int_t^T Z_s^l dB_s-\brak{K_T^l-K_t^l}+\brak{R_T^l-R_t^l},\quad\shoveright{ \quad0\leq t\leq T.} 
\end{equation} 
Indeed, taking conditional expectation on both sides of \eqref{eq;SDGBSDE}, we can obtain
$$\Yb_t^l=V_t^l=\Ebh_t\edg{\xi^l+\int_t^Tf_s^lds},~\mbox{and}~~Y_t^l=V_t^l+\sup_{s\in[t,T]}L_s^l(V_s),~1\leq l\leq N.$$
Then we can check the constraint of the $N$-dimensional reflected $G$-BSDE:
\begin{align*}
	\sum_{l=1}^N\theta^l \Ebh\edg{-Y_t^l}&=\sum_{l=1}^N\theta^l \Ebh\edg{-\brak{V_t^l+\sup_{s\in[t,T]}\Lv^l(V_s)}}\\
	&\leq -\sum_{l=1}^N\theta^l \Lv^l(V_t)+\sum_{l=1}^N\theta^l \Ebh\edg{-V_t^l}=-\brak{\sum_{l=1}^N\theta^l \Ebh\edg{-V_t^l}}^++\sum_{l=1}^N\theta^l \Ebh\edg{-V_t^l}\leq0.
\end{align*}

By the definition \eqref{eq:constR} of $R$, we have
$\sup_{s\in[t,T]}\Lv^l(V_s)=\Lv^l(V_t),~dR_t^l-a.e.,~\mbox{and}~\1_{\crl{\Lv^l(V_t)=0}}=0,~dR_t^l-a.e..$ Recalling \eqref{eq:L-l-t}, we find that the measures $d\abv{R_t}$ and $dR^l_t$ are equivalent for $1\leq l\leq N$ and $\1_{\crl{\Lv^l(V_t)>0}}=\1_{\crl{\abv{\Lv(V_t)}>0}},~d\abv{R_t}-a.e.,$ whenever $\theta^l>0$.
Then we can always find such $j\in\crl{1,\cdots,N}$ by the assumption ($H_\theta$) and deduce that
\begin{align*}
	&\quad\int_0^T\brak{\sum_{l=1}^N\theta^l \Ebh\edg{-Y_t^l}}d\abv{R_t}\\
	&=\int_0^T\brak{\sum_{l=1}^N\theta^l \Ebh\edg{-\brak{V_t^l+\Lv^l(V_t)}}}d\abv{R_t}=\int_0^T\brak{-\sum_{l=1}^N\theta^l \Lv^l(V_t)+\sum_{l=1}^N\theta^l \Ebh\edg{-V_t^l}}\1_{\crl{\Lv^j(V_t)>0}}d\abv{R_t}=0.
\end{align*}

\vspace{0.2em}
\noindent{\bf Step 2. Uniqueness.} 

 The uniqueness is proved by contradiction. When $\theta^l= 0$, the uniqueness of the $l$-th equation is trivial. Otherwise, $\theta^l>0$, $1\leq l\leq N$. Suppose that $\brak{Y^1,Z^1,K^1,R^1}$ and $\brak{Y^2,Z^2,K^2,R^2}$ are two different solutions to the MR$G$-BSDE \eqref{eq:constdriver}. Reciprocally, we can construct two solutions for the $G$-BSDE \eqref{eq;SDGBSDE} index by $1\leq l \leq N$,  $$\brak{Y^{1,l}_\cdot-\brak{R^{1,l}_T-R^{1,l}_\cdot}, Z^{1,l}_\cdot, K^{1,l}_\cdot}\quad \mbox{and}\quad \brak{Y_\cdot^{2,l}-\brak{R_T^{2,l}-R_\cdot^{2,l}},Z^{2,l}_\cdot,K_\cdot^{2,l}}.$$
From Theorem \ref{thm:1deu}, we have on one hand the uniqueness of the solution to one-dimensional $G$-BSDEs \eqref{eq;SDGBSDE}, which implies that the two solutions above identical in appropriate solution spaces. 
In particular, $Y^{1,l}_t-\brak{R^{1,l}_T-R^{1,l}_t} = Y_t^{2,l}-\brak{R_T^{2,l}-R_t^{2,l}}$, $0\leq t\leq T$. 
On the other hand, for $1\leq l\leq N$, since $Y_T^{1,l}=Y_T^{2,l}$, if there exists $t_1<T$ such that $Y_{t_1}^{1,l}\neq Y_{t_1}^{2,l}$, say
$Y_{t_1}^{1,l}> Y_{t_1}^{2,l}$, then $R_T^{1,l}-R_{t_1}^{1,l}>R^{2,l}_T-R^{2,l}_{t_1}.$
Set
$$t_2:=\inf\crl{t>t_1:R_T^{1,l}-R_t^{1,l}=R^{2,l}_T-R^{2,l}_t}.$$
Clearly, $t_2\leq T$ and by \eqref{eq:L-l-t}, $t_2$ does NOT depend on $l$. Thus we have by \eqref{eq:L-l-t} and \eqref{eq:constR} on $t\in\left[t_1,t_2\right) $,
\begin{equation}\label{eq;contra}
\abv{R_T^{1}}-\abv{R_t^{1}}>\abv{R^{2}_T}-\abv{R^{2 }_t}~~\mbox{causing that}~~Y_t^{1,l}>Y_t^{2,l}, ~\mbox{for all}~ 1\leq l\leq N.
\end{equation}
By monotonicity and constraints on expectations, we find
$$\sum_{l=1}^N\theta^l \Ebh\edg{-Y_t^{1,l}}<\sum_{l=1}^N\theta^l \Ebh\edg{-Y_t^{2,l}}\leq 0,~t_1\leq t<t_2,~\mbox{and}~~\sum_{l=1}^N\theta^l \Ebh\edg{-Y_{t_2}^{1,l}}=\sum_{l=1}^N\theta^l \Ebh\edg{-Y_{t_2}^{2,l}}\leq 0,
$$
which implies by flatness of $\brak{Y^1,Z^1,K^1,R^1}$, that
$d\abv{R^1_t}=0, ~\forall t\in \left[t_1,t_2\right).$
Eventually,
$$\abv{R_T^2}-\abv{R_{t_1}^2}\geq \abv{R_T^2}-\abv{R_{t_2}^2}=\abv{R_T^1}-\abv{R_{t_2}^1}=\abv{R_T^1}-\abv{R_{t_1}^1},$$
which is a contradiction to \eqref{eq;contra}.

\vspace{0.2em}
\noindent{\bf Step 3. Check that $Y\in S_G^\a(0,T;\R^N)$ and $Y\in M_G^\b(0,T;\R^N)$.}

For each $1\leq l\leq N$, we have
$$
\abv{Y_t^l}=\abv{\Yb^l_t+\sup_{s\in[t,T]}\Lv^l(V_s)}\leq \abv{\Yb_t^l}+\sup_{s\in[t,T]}\abv{\Lv^l(V_s)},$$
and by \eqref{eq:L-l-t},
\begin{align*}
	\sup_{s\in[t,T]}\abv{\Lv^l(V_s)}
	&\leq \sup_{s\in[t,T]}\frac{\theta^l }{\sum_{j=1}^N(\theta^j)^2}\sum_{j=1}^N\theta^j\Ebh\edg{\abv{V_s^j}}\leq \frac{\theta^l }{\sum_{j=1}^N(\theta^j)^2}\sum_{j=1}^N\theta^j\sup_{s\in[t,T]}\Ebh\edg{\abv{V_s^j}}\\
	&\leq \frac{\theta^l }{\sum_{j=1}^N(\theta^j)^2}\sum_{j=1}^N\theta^j\Ebh\edg{\abv{\xi^j}+\int_0^T\abv{f^j_s}ds}.
\end{align*}
Then standard argument similar to Lemma 3.2 in \cite{LiumdGBSDE} leads to the desired result.
\qed   

\subsection{Case II}
Consider an $N$-dimensional multi-variate MR$G$-BSDE of the following type:
\begin{equation}\label{eq:linearY}
	\left\{\begin{lgathered}
		Y_t^l=\xi^l +\int_t^T f^l(s,Z_s^l)ds-\int_t^T Z_s^l dB_s-\brak{K_T^l-K_t^l}+\brak{R_T^l-R_t^l},
		~1\leq l\leq N; 
		\\
		\sum_{l=1}^N\theta^l \Ebh\edg{-Y_t^l}\leq 0, ~0\leq t\leq T;~\int_0^T \brak{\sum_{l=1}^N\theta^l \Ebh\edg{-Y^l_t}} d\abv{R_t}=0;
	\end{lgathered}	\right.
\end{equation}
where the generators $f^l(t,\om,z^l): [0,T]\x\Om_T\x\R^d\to\R, ~1\leq l\leq N$ satisfy the following assumptions:
\begin{enumerate}
	\item[($H_f^{II}$)] 
	\begin{enumerate}[($i$)]
		\item for each $z^l\in\R^d$, $f^l(\cdot,\cdot,z^l)\in M_G^\b(0,T)$ with $\b>2$ in ($H_\xi$);
		\item there exists some $L>0$ such that, for any $z_1^l,z_2^l\in\R^d$, 
		$\abv{f^l(t,z_1^l)-f^l(t,z_2^l)}\leq L\abv{z_1^l-z_2^l}.$
	\end{enumerate}
\end{enumerate}
\begin{Theorem}\label{theorem:linearY}
	Suppose that ($H_\theta$)-($H_\xi$)-($H_f^{II}$) hold for some $\b>2$ and $L>0$. Then for any $2\leq\a<\b$, MR$G$-BSDE \eqref{eq:linearY}  with parameter $(\xi,f)$ admits a unique solution $(Y,Z,K,R)\in\Sc_G^\a(0,T;\R^N)\x\Ac_D(0,T;\R^N)$. Furthermore, $Y\in M_G^\b(0,T;\R^N)$.
\end{Theorem}

\noindent{\bf \underline{Proof.}}

\vspace{0.2em}
\noindent{\bf Step 1. Existence.}

For each $l\in\crl{1,\cdots,N}$, the standard one-dimensional $G$-BSDE with parameter $(\xi^l,f^l)$:
\begin{equation}\label{eq:constdriverZ}
	\Yt_t^l= \xi^l+\int_t^Tf^l(s,\Zt_s^l)ds-\int_t^T\Zt_s^ldB_s-\brak{\Kt_T^l-\Kt_t^l},~t\in[0,T],
\end{equation}
admits a unique solution $(\Yt^l,\Zt^l,\Kt^l)\in\Sc_G^\a(0,T)$, for any $2\leq\a<\b$ by Theorem \ref{thm:1deu}. 
%{\color{red}Remove conditional expectation.} 
According to Proposition \ref{proposition:LconSG}, $\Lv(\Yt_t)$ is bounded and continuous in $t$. Then putting
$$R_t^l:=\sup_{s\in[0,T]}\Lv^l(\Yt_s)-\sup_{s\in[t,T]}\Lv^l(\Yt_s), ~\mbox{and}~~
	\brak{Y^{l}_t,Z^{l}_t,K^{l}_t}:=\brak{\Yt_t^l+\sup_{s\in[t,T]}\Lv^l(\Yt_s),\Zt_t^l,\Kt_t^l},
$$
one can check that $(Y,Z,K,R)$ is a flat solution.
 
\vspace{0.2em}
\noindent{\bf Step 2. Uniqueness.}

With the help of the uniqueness of the solution $(\Yt^l,\Zt^l,\Kt^l)$ to the $G$-BSDE \eqref{eq:constdriverZ}, for each $1\leq l\leq N$, the uniqueness of $(Y,Z,K,R)$ follows according to an analogous discussion to Step 2 of Theorem \ref{theorem:constdriver}.

\vspace{0.2em}
\noindent{\bf Step 3. Check that $Y\in S_G^\a(0,T;\R^N)$ and $Y\in M_G^\b(0,T;\R^N)$.}

For each $1\leq l\leq N$, we have by \eqref{eq:L-l-t}
$$
\abv{Y_t^l}=\abv{\Yt^l_t+\sup_{s\in[t,T]}\Lv^l(\Yt_s)}\leq \abv{\Yt_t^l}+\sup_{s\in[t,T]}\abv{\Lv^l(\Yt_s)}\leq\abv{\Yt_t^l}+\frac{\theta^l }{\sum_{j=1}^N(\theta^j)^2}\sum_{j=1}^N\theta^j\sup_{s\in[t,T]}\Ebh\edg{\abv{\Yt_s^j}}.$$
Then applying Proposition \ref{prop:1daprY} and following a similar procedure as Step 3 in Theorem \ref{theorem:constdriver}, we can find the desired estimation.
\qed

\subsection{Proof of Theorem \ref{thm:existuniq}}\label{sec:proof-of-theorem-refthmexistuniq}
In this section, we apply a fixed-point technique inspired by Theorem 3.1 of Liu \cite{LiumdGBSDE}. Consider the $N$-dimensional multi-variate MR$G$-BSDE of the following form:
\begin{equation}\label{eq:mdrgb}
	\left\{\begin{lgathered}
		Y_t^l=\xi^l +\int_t^T f^l (s,Y_s,Z_s^l)ds+\-\int_t^T Z_s^ldB_s-\brak{K_T^l-K_t^l}+\brak{R_T^l-R_t^l},~1\leq l\leq N; 
		\\
		\sum_{l=1}^N\theta^l \Ebh\edg{-Y_t^l}\leq 0, ~0\leq t\leq T;~\int_0^T \brak{\sum_{l=1}^N\theta^l \Ebh\edg{-Y^l_t}} d\abv{R_t}=0;
	\end{lgathered}	\right.
\end{equation}
where  the generators
$ f^l(t,\om,y,z^l):[0,T]\times\Om_T\times\R^N\times\R^d\to\R, ~  1\leq l\leq N,$
satisfy the following assumptions:

\begin{enumerate}
	\item[($H_f^*$)] 
	\begin{enumerate}[($i$)]
		\item for each $(y,z^l)\in\R^N\x\R$,  $f^l(\cdot,\cdot,y,z^l)\in M_G^\beta(0,T)$, with $\b>2$ in ($H_\xi$);
		\item there exists some $L>0$ such that, for any $y_1,y_2\in\R^N$ and $z_1^l,z_2^l\in\R^d$,
		
		$\abv{f^l(t,y_1,z_1^l)-f^l(t,y_2,z_2^l)}\leq L\brak{\abv{y_1-y_2}+\abv{z_1^l-z_2^l}}.$
	\end{enumerate}
\end{enumerate}

\begin{Theorem}\label{theorem:hh-contr}
	Suppose that ($H_\theta$)-($H_\xi$)-($H_f^*$) hold for some $\b>2$ and $L>0$. Then there exists a constant  $0<\Dh\leq T$ depending only on $T$, $G$, $L$, $\b$, and $N$ such that for any $h\in\left(0,\Dh\right]$, the $N$-dimensional multi-variate MR$G$-BSDE on the interval $[T-h,T]$:
	\begin{equation}\label{eq:hh-V}
		\left\{\begin{lgathered}
			Y_t^l=\xi^l +\int_t^T f^l(s,Y_s,Z_s^l)ds-\int_t^T Z_s^l dB_s-\brak{K_T^l-K_t^l}+\brak{R_T^l-R_t^l},
			~1\leq l\leq N; 
			\\
			\sum_{l=1}^N\theta^l \Ebh\edg{-Y_t^l}\leq 0, ~T-h\leq t\leq T; ~\int_{T-h}^T \brak{\sum_{l=1}^N\theta^l \Ebh\edg{-Y^l_t}} d\abv{R_t}=0
		\end{lgathered}\right.
	\end{equation}
	admits  a unique solution $(Y,Z,K,R)\in\Sc_G^\a(T-h,T;\R^N)\x\Ac_D(T-h,T;\R^N)$ for each $2\leq\a<\b$. Moreover, $Y\in M_G^\b(T-h,T;\R^N)$.
\end{Theorem}

In order to prove Theorem \ref{theorem:hh-contr}, we need to construct a contraction mapping on $M_G^\b(T-h,T;\R^N)$, where $h$ and $\delta$ will be determined by Lemma \ref{lemma:Hh-def} later. Take $Q\in M_G^\b(T-h,T;\R^N)$ and consider the following MR$G$-BSDE on $[T-h,T]$:
\begin{equation}\label{eq:q-V}
	\left\{\begin{lgathered}
		Y_t^{Q,l}=\xi^l+\int_t^Tf^{Q,l}(s,Z_s^{Q,l})ds-\int_t^TZ_s^{Q,l}dB_s-(K^{Q,l}_T-K^{Q,l}_t)+(R^{Q,l}_T-R_t^{Q,l}),~1\leq l\leq N;\\
		\sum_{l=1}^N\theta^l \Ebh\edg{-Y^{Q,l}_t}\leq 0,~T-h\leq t\leq T;~\int_{T-h}^T \brak{\sum_{l=1}^N\theta^l \Ebh\edg{-Y^{Q,l}_t}} d\abv{R_t^Q}=0.
	\end{lgathered}\right.
\end{equation}
Using Lemma \ref{lemma:funcappr} we have that the generators  $f^{Q,l}(s,z^l):=f^l(s,Q^1_s,\cdots,Q^N_s, z^l), ~1\leq l\leq N$ satisfy condition ($H_f^{II}$). 
Also we denote that $X^Q:=\brak{X^{Q,1},\cdots,X^{Q,N}}^T$, for $X=Y,Z,K$ or $R$.

	By Theorem \ref{theorem:linearY}, we know that equation \eqref{eq:q-V} admits a unique flat solution $(Y^Q,Z^Q,K^Q,R^Q)\in\Sc_G^\a(T-h,T;\R^N)\x\Ac_D(T-h,T;\R^N)$ and $Y^Q\in M_G^\b(T-h,T;\R^N)$. 
Meanwhile, 
$$\brak{Y^{Q,l}_t,Z_t^{Q,l},K_t^{Q,l}}=\brak{\Yb_t^{Q,l}+\brak{R_T^{Q,l}-R^{Q,l}_t},\Zb_t^{Q,l},\Kb_t^{Q,l}},~1\leq l\leq N,$$
where the triples $\crl{(\Yb^{Q,l},\Zb^{Q,l},\Kb^{Q,l})}_{1\leq l\leq N}$ are the unique $\Sc_G^\a(T-h,T)$-solutions to the following $N$ one-dimensional $G$-BSDEs:
$$\Yb_t^{Q,l}=\xi^l+\int_t^T f^{Q,l}(s,\Zb_s^{Q,l})ds-\int_t^T \Zb_s^{Q,l} dB_s-\brak{\Kb^{Q,l}_T-\Kb_t^{Q,l}},~t\in[T-h,T],$$
As Case II, we conclude from the construction of solutions that 
$$
	\Yb^{Q,l}_t=\Ebh_t\edg{\xi^l+\int_t^T f^l(s,Q_s,Z_s^{Q,l})ds}, ~\forall1\leq l\leq N,$$
and
$$
R_t^{Q,l}=\sup_{s\in[T-h,T]}\Lv^l(\Yb^Q_s)-\sup_{s\in[t,T]}\Lv^l(\Yb_s^Q),~\forall t\in[T-h,T].
$$
Then we could define a mapping $\GH:M_G^\b(T-h,T;\R^N)\to M_G^\b(T-h,T;\R^N)$ by
$$\GH(Q):=Y^Q,~ \forall Q\in M_G^\b(T-h,T;\R^N).$$

\begin{Lemma}\label{lemma:Hh-def}
	There exists some $\Dh>0$ such that for each $h\in\left(0,\Dh\right]$, $\GH$ is a contraction under the norm $\norm{\cdot}_{M_G^\b(T-h,T;\R^N)}$.
\end{Lemma}

\noindent{\bf \underline{Proof.}}
Suppose $Q^i\in M_G^\b(T-h,T;\R^N),  ~i=1,2$. For each $i=1,2$, denote by $$\brak{Y^{Q^i},Z^{Q^i},K^{Q^i},R^{Q^i}}\in\Sc_G^\a(T-h,T;\R^N)\x\Ac_D(T-h,T;\R^N)$$ the unique solution to MR$G$-BSDE \eqref{eq:q-V} with parameter $(\xi, f^{Q^i})$.

Since each triple of processes $\crl{(\Yb^{Q^i,l},\Zb^{Q^i,l},\Kb^{Q^i,l})}_{1\leq l\leq N}$ uniquely solves the standard $G$-BSDEs 
$$\Yb_t^{Q^i,l}=\xi^l+\int_t^T f^{Q^i,l}(s,\Zb_s^{Q^i,l})ds-\int_t^T \Zb_s^{Q^i,l} dB_s-\brak{\Kb^{Q^i,l}_T-\Kb_t^{Q^i,l}},~t\in[T-h,T],$$
we have
$$\brak{Y^{Q^i,l}_t,Z_t^{Q^i,l},K_t^{Q^i,l}}=\brak{\Yb_t^{Q^i,l}+\sup_{s\in[t,T]}L_s^l(\Yb_s^{Q^i}),\Zb_t^{Q^i,l},\Kb_t^{Q^i,l}},~1\leq l\leq N.$$
By Proposition \ref{prop:1daprDY}, there exists a constant $C_{(\b)}>0$ depending on $T$, $G$, $\b$ and $L$ such that for each $t\in[0,T]$ ,
$$	\abv{\Yb^{Q^1,l}_t-\Yb^{Q^2,l}_t}^\b\leq C_{(\b)}\Ebh_t\edg{\brak{\int_t^T L\abv{Q^{1}_u-Q^{2}_u}du}^\b}= C_{(\b)} L^\b\Ebh_t\edg{\brak{\int_{T-h}^T \abv{Q^{1}_u-Q^{2}_u}du}^\b}.$$

From Equation \eqref{eq:L-l-t}, we obtain
\begin{align}\label{eq:ytbar}
	&\quad\abv{Y^{Q^1,l}_t-Y^{Q^2,l}_t}\notag\\
	&\leq \abv{\Yb^{Q^1,l}_t-\Yb^{Q^2,l}_t}+\abv{\sup_{s\in[t,T]}\Lv^l(\Yb_s^{Q^1})-\sup_{s\in[t,T]}\Lv^l(\Yb_s^{Q^2})}\leq \abv{\Yb^{Q^1,l}_t-\Yb^{Q^2,l}_t}+\sup_{s\in[t,T]}\abv{\Lv^l(\Yb_s^{Q^1})-\Lv^l(\Yb_s^{Q^2})}\notag\\
	&\leq \abv{\Yb^{Q^1,l}_t-\Yb^{Q^2,l}_t}+\sup_{s\in[t,T]}\frac{\theta^l }{\sum_{j=1}^N(\theta^j)^2}\sum_{j=1}^N\theta^j\Ebh\edg{\abv{\Yb^{Q^1,j}_s-\Yb^{Q^2,j}_s}}\notag\\
	&\leq \abv{\Yb^{Q^1,l}_t-\Yb^{Q^2,l}_t}+\frac{\theta^l }{\sum_{j=1}^N(\theta^j)^2}\sum_{j=1}^N\theta^j\crl{\sup_{s\in[t,T]}\Ebh\edg{\abv{\Yb^{Q^1,j}_s-\Yb^{Q^2,j}_s}}}.
\end{align} 
  
Then for $t\in [T-h, T]$,
\begin{align*}
	&\quad\Ebh\edg{\abv{Y^{Q^1,l}_t-Y^{Q^2,l}_t}^\b}\\
	&\leq 2^{\b-1}\Ebh\edg{\abv{\Yb^{Q^1,l}_t-\Yb^{Q^2,l}_t}^\b}+\frac{2^{\b-1}N^{\b-1}\brak{\theta^l}^\b}{\edg{\sum_{j=1}^N\brak{\theta^j}^2}^\b}\sum_{j=1}^N\brak{\theta^j}^\b\crl{\sup_{s\in[t,T]} \Ebh\edg{\abv{\Yb^{Q^1,j}_s-\Yb^{Q^2,j}_s}^\b}}\\
	&\leq2^{\b-1}C_{(\b)} L^\b\Bigg\{\Ebh\edg{\Ebh_t\edg{\brak{\int_{T-h}^T \abv{Q^{1}_u-Q^{2}_u}du}^\b}}\\
	&~~~~~~~~~~~~~~~~~~~~~~~~~~~~~~~~~~~~~~~~+\frac{N^{\b-1}(\theta^l)^\b}{\edg{\sum_{j=1}^N\brak{\theta^j}^2}^\b}\sum_{j=1}^N\brak{\theta^j}^\b\crl{\sup_{s\in[t,T]} \Ebh\edg{\Ebh_s\edg{\brak{\int_{T-h}^T \abv{Q^{1}_u-Q^{2}_u}du}^\b}}}\Bigg\}\\
	&\leq 2^{\b-1} h^{\b-1}C_{(\b)}L^\b\crl{1+\frac{N^{\b-1}(\theta^l)^\b}{\edg{\sum_{j=1}^N\brak{\theta^j}^2}^\b}\sum_{j=1}^N\brak{\theta^j}^\b}\Ebh\edg{\int_{T-h}^T\abv{Q^{1}_u-Q^{2}_u}^\b du}.
\end{align*}
Here the inequality holds uniformly for $t\in[T-h,T]$ and $1\leq l\leq N$, leading to that
\begin{align*}
	\norm{Y^{Q^1}-Y^{Q^2}}_{M_G^\b(T-h,T;\R^N)}^\b&\leq N^{\frac\b2-1}\int_{T-h}^T\brak{\sum_{l=1}^N\Ebh\edg{\abv{Y^{Q^1,l}_t-Y^{Q^2,l}_t}^\b}} dt\\
	&\leq 2^{\b-1} N^{\frac\b2}h^\b C_{(\b)} L^\b\crl{1+\frac{N^{\b-2}\edg{\sum_{j=1}^N\brak{\theta^j}^\b}^2}{\edg{\sum_{j=1}^N\brak{\theta^j}^2}^\b}} \norm{Q^1-Q^2}_{M_G^\b(T-h,T;\R^N)}^\b.
\end{align*}
Then we can define
$$\Dh:=\crl{2^{\frac1\b-2}N^{-\frac12}C^{-\frac1\b}_{(\b)}L^{-1}\crl{1+\frac{N^{\b-2}\edg{\sum_{j=1}^N\brak{\theta^j}^\b}^2}{\edg{\sum_{j=1}^N\brak{\theta^j}^2}^\b}}^{-\frac1\b}}\wedge T,$$
such that $\norm{Y^{Q^1}-Y^{Q^2}}_{M_G^\b(T-h,T;\R^N)}\leq\frac12\norm{Q^1-Q^2}_{M_G^\b(T-h,T;\R^N)},~\forall h\in\left(0,\Dh\right]$.\qed

\vspace{1em}

\noindent{\bf \underline{Proof of Theorem \ref{theorem:hh-contr} }}

\vspace{0.2em}
\noindent{\bf Step 1. Existence.} 

Take $\Dh$ as in Lemma \ref{lemma:Hh-def} and let $0<h\leq \Dh$.
Put $Y^{[0]}\equiv 0$, $X^{[i+1],l}=X^{Y^{[i]},l}, 1\leq l\leq N$ for $X=Y,Z,K,R$, which are iteratively defined by solving the  $N$-dimensional multi-variate MR$G$-BSDEs:
\begin{equation*}
	\left\{\begin{lgathered}
		Y_t^{[i+1],l}=\xi^l +\int_t^T f^{Y^{[i]},l}(s,Z_s^{[i+1],l})ds-\int_t^T Z_s^{[i+1],l} dB_s-\brak{K^{[i+1],l}_T-K_t^{[i+1],l}}\\
		~~~~~~~~~~~~+\brak{R_T^{[i+1],l}-R_t^{[i+1],l}},
		~1\leq l\leq N; 
		\\
		\sum_{l=1}^N\theta^l \Ebh\edg{-Y_t^{[i+1],l}}\leq 0,~T-h\leq t\leq T;~\int_{T-h}^T \brak{\sum_{l=1}^N\theta^l \Ebh\edg{-Y^{[i+1],l}_t}} d\abv{R_t^{[i+1]}}=0.
	\end{lgathered}\right.
\end{equation*}
From Lemma \ref{lemma:Hh-def}, for each $i\in\N$,
\begin{align*}
	\norm{Y^{[i+1]}-Y^{[i]}}_{M_G^\b(T-h,T;\R^N)}\leq\frac12\norm{Y^{[i]}-Y^{[i-1]}}_{M_G^\b(T-h,T;\R^N)}\leq\cdots\leq \frac{1}{2^i}\norm{Y^{[1]}-Y^{[0]}}_{M_G^\b(T-h,T;\R^N)},
\end{align*}
from which we conclude that $\crl{Y^{[i]}}_{i\geq 0}\subset M_G^\b(T-h,T;\R^N)$ is a Cauchy sequence.

While at the $i$-th iteration , we obtain that 
$$\brak{Y^{[i],l}_t,Z_t^{[i],l},K_t^{[i],l}}=\brak{\Yb_t^{[i],l}+\brak{R^{[i],l}_T-R^{[i],l}_t},\Zb_t^{[i],l},\Kb_t^{[i],l}},~1\leq l\leq N,$$
where $(\Yb^{[i],l},\Zb^{[i],l},\Kb^{[i],l})$ is the unique $\Sc_G^\a(T-h,T)$-solution to the standard $G$-BSDE indexed by $l$:
$$\Yb_t^{[i],l}=\xi^l+\int_t^T f^{Y^{[i-1]},l}(s,\Zb_s^{[i],l})ds-\int_t^T \Zb_s^{[i],l} dB_s-\brak{\Kb^{[i],l}_T-\Kb_t^{[i],l}}.$$

For each $1\leq l\leq N$, estimate \eqref{eq:ytbar} yields that
$$\abv{Y^{[i+1],l}_t-Y^{[j+1],l}_t}\leq\abv{\Yb^{[i+1],l}_t-\Yb^{[j+1],l}_t}+\frac{\theta^l }{\sum_{k=1}^N\brak{\theta^k}^2}\sum_{k=1}^N\theta^k\crl{\sup_{s\in[t,T]}\Ebh\edg{\abv{\Yb^{[i+1],k}_s-\Yb^{[j+1],k}_s}}}.$$
And by Proposition \ref{prop:1daprDY}, there exists some constant $C_{(\a)}>0$ depending on $T$, $G$, $L$ and $\a$ such that
\begin{equation}\label{eq:Dytbar}
	\abv{\Yb^{[i+1],l}_t-\Yb^{[j+1],l}_t}^\a\leq C_{(\a)} L^\a\Ebh_t\edg{\brak{\int_{T-h}^T \abv{Y^{[i]}_s-Y^{[j]}_s}ds}^\a}.
\end{equation}
Then we obtain that the left-hand side of \eqref{eq:Dytbar} is a Cauchy sequence in $\Sc_G^\a(T-h,T;\R^N)$ by showing that
\begin{align*}
	&\quad\abv{Y^{[i+1],l}_t-Y^{[j+1],l}_t}^\a\\
	&\leq 2^{\a-1}\abv{\Yb^{[i+1],l}_t-\Yb^{[j+1],l}_t}^\a+\frac{2^{\a-1}N^{\a-1}(\theta^l)^\a}{\edg{\sum_{k=1}^N\brak{\theta^k}^2}^\a}\sum_{k=1}^N\brak{\theta^k}^\a\crl{\sup_{s\in[t,T]}\Ebh\edg{\abv{\Yb^{[i+1],k}_s-\Yb^{[j+1],k}_s}^\a}}\\
	&\leq 2^{\a-1}C_{(\a)}L^\a\crl{\Ebh_t\edg{\brak{\int_{T-h}^T\abv{Y^{[i]}_s-Y^{[j]}_s}ds}^\a}+ \frac{N^{\a-1}(\theta^l)^\a\sum_{k=1}^N\brak{\theta^k}^\a}{\edg{\sum_{k=1}^N\brak{\theta^k}^2}^\a}\Ebh\edg{\brak{\int_{T-h}^T\abv{Y^{[i]}_s-Y^{[j]}_s}ds}^\a}}.
\end{align*}
Applying Theorem \ref{theorem:simpleDoob} and the H\"older inequality, there exists some constant $C'>0$ depending on $\a$, $\b$ and $G$, such that
\begin{align*}
	&\quad\Ebh\edg{\sup_{t\in[T-h,T]}\abv{Y^{[i+1],l}_t-Y^{[j+1],l}_t}^\a}\\
	&\leq 2^{\a-1}C_{(\a)}L^\a\Bigg\{\Ebh\edg{\sup_{t\in[T-h,T]}\Ebh_t\edg{\brak{\int_{T-h}^T\abv{Y^{[i]}_s-Y^{[j]}_s}ds}^\a}}\\
	&~~~~~~~~~~~~~~~~~~~~~~~~~~~~~~~~~~~~~~~~~~~~~~~~~~~~~+ \frac{N^{\a-1}(\theta^l)^\a\sum_{k=1}^N\brak{\theta^k}^\a}{\edg{\sum_{k=1}^N\brak{\theta^k}^2}^\a}\Ebh\edg{\brak{\int_{T-h}^T\abv{Y^{[i]}_s-Y^{[j]}_s}ds}^\a}\Bigg\}\\
	&\leq2^{\a-1}C_{(\a)}L^\a T^{\a-\frac\a\b}\crl{C'+\frac{N^{\a-1}(\theta^l)^\a\sum_{k=1}^N\brak{\theta^k}^\a}{\edg{\sum_{k=1}^N\brak{\theta^k}^2}^\a}}\crl{ \Ebh\edg{\int_{T-h}^T\abv{Y^{[i]}_s-Y^{[j]}_s}^\b ds}}^\frac\a\b.
\end{align*}
Then we have 
\begin{align*}
	&\norm{Y^{[i+1]}-Y^{[j+1]}}_{S_G^\a(T-h,T;\R^N)}^\a\leq N^{\ah-1}\sum_{l=1}^N\Ebh\edg{\sup_{t\in[T-h,T]}\abv{Y^{[i+1],l}_t-Y^{[j+1],l}_t}^\a}\\
	&~~~~~~~~~~\leq2^{\a-1}N^{\ah}C_{(\a)}L^\a T^{\a-\frac\a\b}\crl{C'+\frac{N^{\a-2}\edg{\sum_{k=1}^N\brak{\theta^k}^\a}^2}{\edg{\sum_{k=1}^N\brak{\theta^k}^2}^\a}}\norm{Y^{[i]}-Y^{[j]}}_{M_G^\b(T-h,T;\R^N)}^\a.
\end{align*}
Hence $\crl{Y^{[i]}}_{i\geq 0}$ converges to some $Y\in S_G^\a(T-h,T;\R^N)$. Moreover, there is some constant $C_{T, G, N, L, \alpha, \beta}>0$ such that
\begin{equation}\label{eq:estimateYt}
	\norm{Y^{[i]}}_{S_G^\a(T-h,T;\R^N)}\leq C_{T, G, N, L, \alpha, \beta}, ~\mbox{uniformly for all}~i\in\N.
\end{equation}

By the estimate \eqref{eq:Dytbar} and Theorem \ref{theorem:simpleDoob}, we deduce that
\begin{align*}
	&\quad\norm{\Yb^{[i+1]}-\Yb^{[j+1]}}_{S_G^\a(T-h,T;\R^N)}^\a\\
	&\leq N^{\ah-1}\sum_{l=1}^N\Ebh\edg{\sup_{T-h\leq t\leq T}\abv{\Yb^{[i+1],l}_t-\Yb^{[j+1],l}_t}^\a}
	\leq C'C_{(\a)} L^\a N^\ah T^{\a-\frac\a\b}\norm{Y^{[i]}-Y^{[j]}}_{M_G^\b(T-h,T;\R^N)}^\a.
\end{align*}
Consequently, $\crl{\Yb^{[i]}}_{i\geq 0}$ is a Cauchy sequence in $S_G^\a(T-h,T;\R^N)$, thus for some constant $\Cb_{T, G, N, L, \alpha, \beta}>0$,  
\begin{equation}\label{eq:estimateYbt}
	\norm{\Yb^{[i]}}_{S_G^\a(T-h,T;\R^N)}\leq \Cb_{T, G, N, L, \alpha, \beta}, ~\mbox{uniformly for all}~i\in\N.
\end{equation}

Now concerning the convergence of $\crl{Z^{[i]}}_{i\geq0}$, apply Proposition \ref{prop:1daprDZ} and we have
\begin{align*}
	&\quad\Ebh\edg{\brak{\int_{T-h}^T\abv{Z_s^{[i+1],l}-Z_s^{[j+1],l}}^2ds}^\ah}\\
	&\leq C\Bigg\{\norm{\Yb^{[i+1],l}-\Yb^{[j+1],l}}_{S_G^\a(T-h,T)}^\a+\norm{\Yb^{[i+1],l}-\Yb^{[j+1],l}}_{S_G^\a(T-h,T)}^\ah\x\bigg[\norm{\Yb^{[i+1],l}}_{S_G^\a(T-h,T)}^\ah\\
	&\quad\quad\quad+\norm{\Yb^{[j+1],l}}_{S_G^\a(T-h,T)}^\ah+\norm{\int_{T-h}^T\abv{f^{Y^{[i]},l}(s,0)}ds}_{L_G^\a}^\ah+\norm{\int_{T-h}^T\abv{f^{Y^{[j]},l}(s,0)}ds}_{L_G^\a}^\ah\bigg]\Bigg\}\\
	&\leq C\Bigg\{\norm{\Yb^{[i+1],l}-\Yb^{[j+1],l}}_{S_G^\a(T-h,T)}^\a+\norm{\Yb^{[i+1],l}-\Yb^{[j+1],l}}_{S_G^\a(T-h,T)}^\ah\x\bigg[\norm{\Yb^{[i+1],l}}_{S_G^\a(T-h,T)}^\ah\\
	&\quad\quad\quad+\norm{\Yb^{[j+1],l}}_{S_G^\a(T-h,T)}^\ah+\norm{Y^{[i]}}_{S_G^\a(T-h,T;\R^N)}^\ah+\norm{Y^{[j]}}_{S_G^\a(T-h,T;\R^N)}^\ah+\norm{\int_{T-h}^T\abv{f^{l}(s,\0,0)}ds}_{L_G^\a}^\ah\bigg]\Bigg\}\\
	&\leq C\crl{\norm{\Yb^{[i+1],l}-\Yb^{[j+1],l}}_{S_G^\a(T-h,T)}^\a+\norm{\Yb^{[i+1],l}-\Yb^{[j+1],l}}_{S_G^\a(T-h,T)}^\ah},
\end{align*}
in which the constant $C$ may vary from line to line and the last inequality is given by assumption ($H_f^*$) and the estimates \eqref{eq:estimateYt}-\eqref{eq:estimateYbt}. 
Therefore $\crl{Z^{[i]}}_{i\geq 0}$ is a Cauchy sequence in $H_G^\a(T-h,T;\R^{N\x d})$,  converging to some $Z\in H_G^\a(T-h,T;\R^{N\x d})$.

Next consider the terms $R^{[i+1],l}-R^{[j+1],l}$, by the constructions we have
 $R^{[i+1],l}_{T-h}=R^{[j+1],l}_{T-h}=0$ and
\begin{align*}
	&\quad\sup_{t\in[T-h,T]}\abv{R^{[i+1],l}_t-R^{[j+1],l}_t}\\
	&\leq \sup_{t\in[T-h,T]}\brak{\abv{\sup_{s\in[T-h,T]}\Lv^l(\Yb^{[i+1]}_s)-\sup_{s\in[T-h,T]}\Lv^l(\Yb^{[j+1]}_s)}+\abv{\sup_{s\in[t,T]}\Lv^l(\Yb^{[i+1]}_s)-\sup_{s\in[t,T]}\Lv^l(\Yb^{[j+1]}_s)}}\\
	&\leq 2\sup_{s\in[T-h,T]}\abv{\Lv_s^l(\Yb^{[i+1]}_s)-\Lv_s^l(\Yb^{[j+1]}_s)}\leq \frac{2\theta^l }{\sum_{k=1}^N\brak{\theta^k}^2}\sum_{k=1}^N\theta^k\norm{\Yb^{[i+1],k}-\Yb^{[j+1],k}}_{S_G^\a(T-h,T)}.
\end{align*}
Therefore, $R^{[i],l}$ converges to some bounded, deterministic, continuous, and non-decreasing function $R^l$ in $t$ with $R^l_{T-h}=0$.
Hereafter, we define the process $K$ by
\begin{align}\label{eq:defKl}
	K^l_t:=Y^l_t+\int_{T-h}^t f^l(s,Y_s,Z_s^l)ds-\int_{T-h}^tZ_s^ldB_s-Y^l_{T-h}+R^l_t, ~T-h\leq t\leq T.
\end{align}
Then the assertion that $K^{[i]}_t\to K_t\in L_G^\a(\Om_t;\R^N)$ and $K\in\Ac_G^\a(T-h,T;\R^N)$ can be acquired through quasi-continuity property of the right-hand side of \eqref{eq:defKl} along with some approximations applying Lemma \ref{lemma:BDG} and convergence results already derived: for any $1\leq l\leq N$ and $i\in\N$,
\begin{align*}
	&\quad\sup_{t\in[T-h,T]}\Ebh\edg{\abv{K^{[i],l}_t-K^l_t}^\a}\\
	&\leq C_{T,L,\a} \crl{\norm{Y^{[i],l}-Y^l}^\a_{S_G^\a(T-h,T)}+\norm{Z^{[i],l}-Z^l}^\a_{H_G^\a(T-h,T)}+\sup_{t\in[T-h,T]}\abv{R^{[i],l}_t-R^l_t}^\a}\\
	&\leq C_{T,L,N,\a,\b,\theta} \brak{\norm{Y^{[i]}-Y}^\a_{S_G^\a(T-h,T;\R^N)}+\norm{Z^{[i],l}-Z^l}^\a_{H_G^\a(T-h,T)}}.
\end{align*}

\vspace{0.2em}
\noindent{\bf Step 2. Uniqueness.}

Let $\brak{Y',Z',K',R'}$ be another solution to the MR$G$-BSDE \eqref{eq:hh-V} on $[T-h,T]$. From Lemma \ref{lemma:Hh-def}, we have $Y=Y'$ because they are both fixed-point solutions to $\Gamma(Y) = Y$ due to the contraction, which leads to the uniqueness.

\begin{Remark}
	Note that the modified generators for obtaining a fixed point in \cite{LiumdGBSDE} take the form (we put the dimension therein as $N$ to avoid notational confusions):
	\begin{equation*}
		f^{l,U}(t,\om,y^l,z^l):=f^l(t,\om,U^1_t,\cdots,U^{l-1}_t,y^l,U^{l+1}_t,\cdots,U^N_t,z^l):~[0,T]\x\Om_T\x\R\x\R^d\to\R,~1\leq l\leq N,
	\end{equation*} 
where $U\in M_G^\b(T-h,T;\R^N)$, of which $Y^{(i)},~i\in\N$ takes the role at the $(i+1)$-th iteration to solve the following $N$ one-dimensional standard $G$-BSDEs:
\begin{equation*}
	Y^{(i+1),l}_t=\xi^l+\int_t^Tf^{l,Y^{(i)}}(s,Y^{(i+1),l}_s,Z^{(i+1),l}_s)ds-\int_t^TZ^{(i+1),l}dB_s-(K^{(i+1),l}_T-K^{(i+1),l}_t),~t\in[T-h,T].
\end{equation*}
 Considering the contraction procedure here, we set generators $f^{Q,l}(t,\om,z^l):=f^l(t,\om,Q^1,\cdots,Q^N,z^l):[0,T]\x\Om_T\x\R^d\to\R,~1\leq l\leq N$ with $Q\in M_G^\b(T-h,T;\R^N)$, which are only diagonally dependent on $z\in\R^{N\x d}$, so as to utilize at each iteration the well-posedness of $N$-dimensional multi-variate MR$G$-BSDE in Case II.  
\end{Remark}
\begin{Theorem}\label{theorem:0Tsolution}
	Suppose that ($H_\theta$)-($H_\xi$)-($H_f^*$) hold for some $\b>2$ and $L>0$. Then for any $~2\leq\a<\beta$, MR$G$-BSDE \eqref{eq:mdrgb} admits a unique solution $(Y,Z,K,R)\in\Sc_G^\a(0,T;\R^N)\x\Ac_D(0,T;\R^N)$. Moreover, $Y\in M_G^\beta(0,T;\R^N)$.
\end{Theorem}
\noindent{\bf \underline{Proof.}}
First we prove the existence. Choose $n\in\N*$ such that $n\delta>T$ and let $h:=\frac T n$. By Theorem \ref{theorem:hh-contr}, the MR$G$-BSDE \eqref{eq:hh-V} admits a unique solution $(Y^n,Z^n,K^n,R^n)\in\Sc_G^\a(T-h,T;\R^N)\x\Ac_D(T-h,T;\R^N)$ and $Y^n\in M_G^\b(T-h,T;\R^N)$. By the quasi-continuity of $Y^{n}_{T-h}$ and the boundedness in $L^\beta_G(\Omega_{T-h})$ given by Proposition \ref{prop:1daprY}, we see that $Y^n_{T-h}\in L^\beta_G(\Omega_{T-h})$. Then we put $Y^{n}_{T-h}$ as the terminal value and apply Theorem \ref{theorem:hh-contr} on the time interval $[T-2h,T-h]$. In this way, the MR$G$-BSDE \eqref{eq:hh-V} with parameter $(Y_{T-h}^{n},f)$ on $T-2h\leq t\leq T-h$ admits a unique solution $(Y^{n-1},Z^{n-1},K^{n-1},R^{n-1})\in\Sc_G^\a(T-2h,T-h;\R^N)\x\Ac_D(T-2h,T-h;\R^N)$. Repeating this procedure until $t\in[0,h]$, we obtain a sequence $\crl{\brak{Y^i,Z^i,K^i,R^i}}_{1\leq i\leq n}$.

Set for each $1\leq l\leq N$,
\begin{flalign*}
	Y_t^l&=\sum_{i=1}^n Y^{i,l}_t\1_{\left[(i-1)h,ih\right)}(t)+\xi^{l}\1_{\crl{T}}(t),~Z_t^l=\sum_{i=1}^n Z^{i,l}_t\1_{\left[(i-1)h,ih\right)}(t)+Z_T^{n,l}\1_{\crl{T}}(t),\\
	K^l_t&=\sum_{i=1}^n K^{i,l}_{ih\wedge t}\1_{\left((i-1)h,T\right]}(t),~~~~~~~~~~~~~~R^l_t=\sum_{i=1}^n R^{i,l}_{ih\wedge t}\1_{\left((i-1)h,T\right]}(t).
\end{flalign*}
One can easily check that $(Y,Z,K,R)\in\Sc_G^\a(0,T;\R^N)\x\Ac_D(0,T;\R^N)$ and is a solution to the MR$G$-BSDE \eqref{eq:mdrgb}. By the same reasoning as in the proof of Theorem \ref{theorem:constdriver}, $Y\in M_G^\b(0,T;\R^N)$.

The uniqueness follows by that on each small time interval $[(i-1)h,ih],1\leq i\leq n$.
\qed
\begin{Remark}
	In this paper the reflections are imposed on convex combinations of $G$-expectations of all $Y$-components $\sum_{l=1}^N\theta^l\Ebh\edg{-Y^l}$ rather than $G$-expectation of convex combinations $\Ebh\edg{-\sum_{l=1}^N\theta^l Y^l}$. In that case, to prove the solvability, we shall encounter a combination of non-increasing asymmetric $G$-martingales, and the property of such process under $G$-expectation is not clear, which is postponed to future research.        
\end{Remark}

\section{Extension: mean reflection related to nonlinear expectation}

In this section, let us extend the results of MR$G$-BSDE under nonlinear expectation framework:
\begin{equation}\label{eq:GBSDEDGEt}
	\left\{\begin{lgathered}
		Y_t^l=\xi^l +\int_t^T f^l (s,Y_s,Z_s^l)ds-\int_t^T Z_s^ldB_s-\brak{K_T^l-K_t^l}+\brak{R_T^l-R_t^l},	~1\leq l\leq N; 
		\\
		\sum_{l=1}^N\theta^l \Et\edg{-Y_t^l}\leq 0,~0\leq t\leq T;~\int_0^T\brak{\sum_{l=1}^N\theta^l \Et\edg{-Y^l_t}} d\abv{R_t}=0;
	\end{lgathered}\right.	
\end{equation}
where $\Et\edg{\cdot}$ is a nonlinear expectation dominated by $G$-expectation $\Ebh\edg{\cdot}$, i.e.,
\begin{equation}\label{eq:dominance}
	\Et\edg{X}-\Et\edg{Y}\leq\Ebh\edg{X-Y},~\forall X,Y\in L_G^1(\Om_T).
\end{equation}
Obviously, we get that for any $X,Y\in L_G^1(\Om_T)$,
\begin{align}
	-\Ebh\edg{-X}\leq(-\Et\edg{-X})\wedge\Et\edg{X}&\leq(-\Et\edg{-X})\vee\Et\edg{X}\leq\Ebh\edg{X},\notag\\
	\mbox{and}~~\abv{\Et\edg{X}-\Et\edg{Y}}&\leq\Ebh\edg{\abv{X-Y}}.\label{eq:Etabv}
\end{align}
Here we define an operator $\Lt:L_G^1(\Om_T;\R^N)\to\R^N$ similarly as $\Lv$ by
$$\Lt(\eta):=\arg\min\crl{\abv{\xb}:-\sum_{j=1}^N\theta^l \xb^j+\sum_{j=1}^N\theta^j\Et\edg{-\eta^j}\leq0}=\frac{\sum_{l=1}^N\theta^l \e_l}{\sum_{j=1}^N\brak{\theta^j}^2}\brak{\sum_{j=1}^N\theta^j\Et\edg{-\eta^j}}^+;$$
and 
$$\Ht(\xb,\eta):=-\sum_{l=1}^N\theta^l \xb^l+\sum_{l=1}^N\theta^l \Et\edg{-\eta^l},~\eta\in L_G^1(\Om_T;\R^N).$$
Following Lemmas 4.1 and 4.2 in Liu and Wang \cite{LWMFGB}, we can obtain that $\Lt$ is well-defined. And Propositions \ref{prop:CEcont} and \ref{proposition:LconSG} still hold with $\Lv$ and $H(\cdot,\cdot)$ replaced by $\Lt$ and $\Ht(\cdot,\cdot)$. Thus we have  
\begin{Theorem}Suppose that ($H_\theta$) hold. We assume that $\xi\in L_G^\b(\Om_T;\R^N)$ and $\sum_{l=1}^N\theta^l\Et\edg{-\xi^l}\leq0$.
	\begin{enumerate}[($i$)]
		\item Let ($H_f^{II}$) hold. Then for any $2\leq\a<\b$, MR$G$-BSDE \eqref{eq:linearY} with $\Et$-reflection and Skorokhod condition of \eqref{eq:GBSDEDGEt} admits a unique solution $(Y,Z,K,R)\in\Sc_G^\a(0,T; \R^N)\x\Ac_D(0,T;\R^N)$. Moreover, $Y\in M_G^\b(0,T;\R^N)$.
		\item Let ($H_f^*$) hold. Then for each $2\leq\a<\b$, MR$G$-BSDE \eqref{eq:GBSDEDGEt}
		admits  a unique solution $(Y,Z,K,R)\in\Sc_G^\a(0,T; \R^N)\x\Ac_D(0,T;\R^N)$.
	\end{enumerate}
\end{Theorem}
%\noindent{\bf\underline{Proof.}}
The proof of the theorem follows the procedures in Theorem \ref{theorem:linearY}, Theorem \ref{theorem:hh-contr}, Lemma \ref{lemma:Hh-def} and Theorem \ref{theorem:0Tsolution}. Note that the operator $\Lt$ shares the same properties as $\Lv$, which defines $R$-terms under the $G$-expectation mean reflection. Moreover, the estimates of the solution under $\Ebh$-norms hold as well due to relations \eqref{eq:Etabv}.

\section*{Acknowledgements}

The authors are grateful to Zihao Gu and Kun Xu for helpful discussions. This work is partially supported by
NSFC under the project No. 12371473 and  by the Tianyuan Fund for Mathematics of NSFC under the project No. 12326603.

\end{document}